# Instability and non-uniqueness in the Cauchy problem for the Euler equations of an ideal incompressible fluid. Part II

Misha Vishik

Department of Mathematics

The University of Texas at Austin

Austin, TX 78712, USA

## § 0. Introduction

In this paper we prove instability of a certain class of radially symmetric flows of an ideal incompressible fluid in plane. In Part I of the present paper we constructed a radially nonsymmetric solution of the Euler equations of an ideal incompressible fluid in $\mathbb{R}^2$ with radially symmetric initial Cauchy data with vorticity bounded in time in Lebesgue space $L^Q(\mathbb{R}^2)$. The external force in vorticity formulation of the equations is radially symmetric in $L^1_{loc}([0,\infty); L^Q(\mathbb{R}^2))$.

The strategy of the construction in Part I was to use a family of linearly instable radially symmetric solutions of a special type to produce a radially nonsymmetric flow by using the fundamental scaling properties of the Euler equations. These solutions to the nonlinear problem are self-similar up to "small" and controlled error terms.

Technically, this construction relies on a linear stability result on existence of radially symmetric flows with the spectrum of small oscillations around it (for, example, in $L^2(\mathbb{R}^2)$) having special features used in our nonlinear construction.

In this paper we give a proof of the linear instability result used in the nonlinear construction in Part I of the paper.

The literature on linear stability / instability of incompressible steady flows of an ideal fluid is vast and we will make no attempt to give a complete account in this introduction. In addition to the sheer volume of the literature on linear stability / instability some of it comes in a state that can hardly be considered proved in mathematical sense. Curiously, mathematical proofs of instability due to the discrete spectrum of the perturbation problem often require a new method for every new case considered.

We mention classical texts by H.Lamb [LA], C.C.Lin [LI] and S.Chandrasekhar [CH].

In his remarkable paper [FA] L.D.Faddeev studied the stability problem for a shear flow in a 2D chanel. He proved that for a monotone velocity profile with an inflexion point the instable spectrum can be "generated" in a complex neighborhood of the inflexion point. This work is related to the work of L.D.Faddeev on scattering theory and on K.O.Friedrichs model in perturbation theory of continuous spectra. See also [F1].

To formulate the main result of the present paper we denote by $R(s)$ the distribution of an angular velocity of a radially symmetric incompressible flow in plane as a function of a radius $s$, and by $G(s)$ the distribution of vorticity as a function of a radius. Let an integer $m \neq 0$ denote an azimuthal wave number.

The unstable discrete spectrum of infinitesimal oscillations about the equilibrium flow is described by the following problem in terms of a stream function of a perturbation $e^{im\theta}\psi(s)$, $\theta$ being a polar angle,

$$-\frac{d^2}{ds^2}\psi - \frac{1}{s}\frac{d}{ds}\psi + \frac{m^2}{s^2}\psi + \frac{B(s)}{R(s)-\mu}\psi = 0 \quad . \tag{0.1}$$

Here $\mu$ is a spectral parameter and (linear) instability of the flow corresponds to the condition $Im\,\mu > 0$.

In this paper we construct a class of vorticity profiles $G(s)$ for which there exists an integer $m \geq 2$ so that problem (0.1) has one and exactly one $\mu$ with $Im\,\mu > 0$ corresponding to a nontrivial perturbation with vorticity in $L^2(\mathbb{R}^2)$.

All the problems (0.1) with $m$ being replaced by $m+1, m+2, m+3, \ldots,$

do not have a nontrivial solutions for a spectral parameter $\mu$ with $Im\,\mu > 0$. The choice of a function space is unimportant for the eigenfunctions corresponding to the discrete spectrum, the corresponding eigenfunction of vorticity perturbation belongs to $L^Q(\mathbb{R}^2)$ for any $Q \in [1, \infty]$. It also belongs to any Holder space $L^Q(\mathbb{R}^2)$, $s \in \mathbb{R}$.

The vorticity profile $G(s)$ decays at infinity like $s^{-\alpha}$, $\alpha \in (0,2)$. This feature is a useful one for applying the rescaling in Part I of the paper.

We get rather complete information about solutions to the problem (0.1) for a special class of profiles by studying possible solutions as functions of a real parameter $m$ as $m$ decreases from infinity to 1. Only integer values of $m$ correspond to a linearized stability problem for the perturbations of vorticity.

# § 1. Linearized stability problem and the eigenfunctions.

## Preliminaries and the formulation of the main theorem

We change the notation as compared to Part I of the paper.

Let $G(|x|)$, $x \in \mathbb{R}^2$, $x = (x_1, x_2)$ be radially symmetric $C^\infty(\mathbb{R}^2)$ distribution of vorticity in plane. Define

$$V(x) = \frac{x^\perp}{|x|^2} \int_0^{|x|} s\, G(s)\, ds, \qquad x^\perp = (-x_2, x_1). \tag{1.1}$$

Then $curl\, u(x) = \partial_1 V_2 - \partial_2 V_1 = G(|x|)$.

We choose $\alpha \in (0,2)$ and assume

$$G(s) = s^{-\alpha},\, for\, s \geq M > 0, \tag{1.2}$$

$M$ being a sufficiently large constant.

Any velocity field (1.1) does (formally) satisfy the Euler equations of an ideal incompressible fluid.

To study the spectrum of small perturbations to this steady-state solution we introduce infinitesimal perturbation

$$w(x) = \partial_x^\perp (e^{im\theta} \psi(|x|)), \quad \partial_x^\perp = (-\partial_2, \partial_1) ; \tag{1.3}$$

$$x_1 = |x|\cos\theta, \quad x_2 = |x|\sin\theta, \quad \theta \in \mathbb{R}/2\pi\mathbb{Z}, \; x \in \mathbb{R}^2, m \in \mathbb{Z}\setminus\{0\}.$$

Let also

$$e^{im\theta} g(|x|) = curl\, w(x) = \partial_1 w_2 - \partial_2 w_1. \tag{1.4}$$

From (1.3), (1.4),

$$e^{im\theta} g(|x|) = \Delta(e^{im\theta} \psi(|x|));$$

$$g(s) = \frac{d^2}{ds^2}\psi + \frac{1}{s}\frac{d}{ds}\psi - \frac{m^2}{s^2}\psi \tag{1.5}$$

We have,

$$\psi(s) = -\frac{1}{2|m|} s^{|m|} \int_s^\infty g(\tau)\tau^{1-|m|}\, d\tau - \frac{1}{2|m|} s^{-|m|} \int_0^s g(\tau)\tau^{1+|m|}\, d\tau. \tag{1.6}$$

Indeed,

$$\frac{d^2}{ds^2}\psi + \frac{1}{s}\frac{d}{ds}\psi - \frac{m^2}{s^2}\psi = -\frac{1}{|m|}|m|s^{|m|-2}\left(-s^{2-|m|}g(s)\right) - \frac{1}{|m|}(-|m|)s^{-|m|-2}\left(s^{2+|m|}g(s)\right)$$

$$+ \frac{1}{2|m|}s^{|m|}s^{-1}\frac{d}{ds}\left(s^{2-|m|}g(s)\right) - \frac{1}{2|m|}s^{-|m|}s^{-1}\frac{d}{ds}\left(s^{2+|m|}g(s)\right)$$

$$= 2g(s) + \left(\frac{2-|m|}{2|m|} - \frac{2+|m|}{2|m|}\right)g(s) = g(s).$$

From (1.4) - (1.6),

$$w(x) = \partial_x^\perp \left( e^{im\theta} \psi(|x|) \right) = \psi'(|x|)e^{im\theta} \frac{x^\perp}{|x|} - im\, \psi(|x|)\, e^{im\theta} \frac{x}{|x|^2} =$$

$$\left( -\frac{|m|}{2|m|} |x|^{|m|-1} \int_{|x|}^\infty g(\tau)\tau^{1-|m|}\, d\tau + \frac{|m|}{2|m|} |x|^{-|m|-1} \int_0^{|x|} g(\tau)\tau^{1+|m|}\, d\tau \right) e^{im\theta} \frac{x^\perp}{|x|}$$

$$- im\, \psi(|x|)\, e^{im\theta} \frac{x}{|x|^2} =$$

$$\left( -\frac{1}{2}|x|^{|m|-2} \int_{|x|}^\infty g(\tau)\tau^{1-|m|}\, d\tau\, e^{im\theta} + \frac{1}{2}|x|^{-|m|-2} \int_0^{|x|} g(\tau)\tau^{1+|m|}\, d\tau\, e^{im\theta} \right) x^\perp$$

$$- im\, \psi\, e^{im\theta} \frac{x}{|x|^2}$$

$$= |x|^{|m|-2} \int_{|x|}^\infty g(\tau)\tau^{1-|m|}\, d\tau\, e^{im\theta} \left( -\frac{1}{2}x^\perp + \frac{i}{2} \text{sign}(m) x \right)$$

$$+ |x|^{-|m|-2} \int_0^{|x|} g(\tau)\tau^{1+|m|}\, d\tau\, e^{im\theta} \left( \frac{1}{2}x^\perp + \frac{i}{2} \text{sign}(m) x \right). \tag{1.7}$$

The spectral problem for small oscillations around the steady state solution $V(x)$ looks as follows:

$$\begin{cases} -(V(x), \partial_x)w - (w, \partial_x)V(x) - \partial_x p = \lambda\, w(x) \\ \text{div}\, w = 0. \end{cases} \tag{1.7'}$$

We are looking for a nontrivial solution $w \neq 0$, $w \in L^2(\mathbb{R}^2)$; $\text{Re}\, \lambda \neq 0$,

where $w$ is satisfying

$$w(R_\theta\, x) = e^{im\theta}\, R_\theta\, w(x), \text{ where } R_\theta = \begin{pmatrix} \cos\theta & -\sin\theta \\ \sin\theta & \cos\theta \end{pmatrix}, \forall\, \theta \in \mathbb{R}/2\pi\mathbb{Z}.$$

Then in $S'(\mathbb{R}^2)$ and in $W^{-1,2}(\mathbb{R}^2)$

$$\text{curl}\, w\, (R_\theta\, x) = e^{im\theta}\, \text{curl}\, w(x).$$

Therefore,

$$(x^\perp, \partial_x)\text{curl}\, w = im\, \text{curl}\, w$$

in $S'(\mathbb{R}^2)$.

Differentiating the equation for eigenfunctions, we get

$$-(V(x), \partial_x)\text{curl}\, w - (w, \partial_x)\text{curl}\, V(x) = \lambda\, \text{curl}\, w.$$

But $V(x) = R(|x|)x^\perp$, where $R(|x|) \in C^\infty(\mathbb{R}^2) \cap S_{1,0}^{-\alpha}(\mathbb{R}^2)$. Therefore, in $S'(\mathbb{R}^2)$

$$(-im\, R(|x|) - \lambda)\, curl\, w = (w, \partial_x)\, G(|x|),$$

where

$$R(s) = \frac{1}{s^2} \int_0^s \tau\, G(\tau) d\tau, \quad s \geq 0. \tag{1.8}$$

But $\{x \to (-im\, R(|x|) - \lambda)\}$ is an invertible multiplier in $S(\mathbb{R}^2)$ and in $S'(\mathbb{R}^2)$. The right side

$$(w, \partial_x)\, G(|x|) \in L^1(\mathbb{R}^2) \cap L^2(\mathbb{R}^2).$$

Therefore,

$$curl\, w(x) = e^{im\theta} g(|x|) \in L^1(\mathbb{R}^2) \cap L^2(\mathbb{R}^2).$$

By Sobolev embedding theorem, $w \in L^Q(\mathbb{R}^2)$, $\forall Q \in (2, \infty)$.

Therefore, a unique stream function $e^{im\theta}\psi(|x|)$, such that $w(x) = \partial_x^\perp(e^{im\theta}\psi(|x|))$, $\psi(0) = 0$,

satisfies $\psi(|x|) \in \dot{C}^{1-\frac{2}{Q}}$, $\forall Q \in (2, \infty)$. In particular,

$$|\psi(|x|)| \leq C|x|^{1-\frac{2}{Q}}, \quad \forall x \in \mathbb{R}^2.$$

The spectral problem for small oscillations around the steady state solution $V(x)$ in terms of vorticity perturbations looks as follows:

$$-im\, R(|x|) g(|x|) + im\, \frac{1}{|x|} G'(|x|) \psi(|x|) = \lambda\, g(|x|). \tag{1.9}$$

Let

$$B(s) = \frac{1}{s} G'(s). \tag{1.10}$$

From (1.8) - (1.10)

$$R(s)g(s) - B(s)\psi(s) = \mu g(s), \tag{1.11}$$

$$\mu = i\frac{\lambda}{m}. \tag{1.12}$$

From (1.2), (1.9), (1.10)

$$B(s) = -\alpha s^{-\alpha-2}, \quad s \geq M; \tag{1.13}$$

$$R(s) = C\alpha s^{-2} + \frac{1}{2-\alpha} s^{-\alpha}, \quad s \geq M \tag{1.14}$$

for an appropriate constant $C$.

For the solutions of (1.7'), (1.11) with $Im\, \mu \neq 0$, $w(x) \in L^2(\mathbb{R}^2)$,

it follows, that

$$\partial_x \left( e^{im\theta} \psi(|x|) \right) \in L^2(\mathbb{R}^2). \tag{1.15}$$

Therefore,

$$\int_0^\infty s^{-1} |\psi(s)|^2 ds < \infty, \quad \int_0^\infty s \left| \frac{d\psi}{ds} \right|^2 ds < \infty,$$

and, consequently, $\psi \in L^\infty(\mathbb{R}^2)$. To see this, we have, for any $s > 0$

$$|\psi(s)^2| = 2 \left| \int_s^\infty s^{-\frac{1}{2}} \psi(s) s^{\frac{1}{2}} \frac{d\psi}{ds} ds \right| \leq C < \infty.$$

From (1.15), for any eigenfunction $w(x) \in L^2(\mathbb{R}^2)$ as in (1.7') with $Im\, \mu \neq 0$

$$|g(s)| \leq C(1+s)^{-\alpha-2}, \quad \forall s \geq 0. \tag{1.16}$$

Using (1.6), (1.16) we obtain

$$|g(s)| \leq C(1+s)^{-2\alpha-2}, \quad \forall s \geq 0.$$

Continuing bootstrap, we arrive at the inequality

$$|g(s)| \leq C(1+s)^{-(|m|+2+\alpha)}, \quad \forall s \geq 0. \tag{1.17}$$

Using (1.6), (1.17), (1.11), (1.13), (1.14), we get

$$\mu g(s) = \alpha s^{-\alpha-2} \left( -\frac{1}{2|m|} s^{-|m|} \left[ \int_0^\infty g(\tau) \tau^{1+|m|} d\tau + O(s^{-\alpha}) \right] + O\left(s^{-|m|-\alpha}\right) \right) + O(s^{-|m|-2-2\alpha})$$

Assuming $\int_0^\infty g(\tau) \tau^{1+|m|} d\tau \neq 0$, this implies

$$\mu = -\frac{\alpha}{2|m|} \int_0^\infty g(\tau) \tau^{1+|m|} d\tau; \tag{1.18}$$

$$\psi(s) = -\frac{1}{2|m|} s^{-|m|} \int_0^\infty g(\tau) \tau^{1+|m|} d\tau + O(s^{-|m|-\alpha}) \text{ as } s \to \infty, \tag{1.19}$$

provided $g(s)$ is normalized as follows:

$$g(s) \sim s^{-|m|-\alpha-2} \quad \text{as } s \to \infty. \tag{1.20}$$

Lemma 1.1. If (1.7'), (1.11) holds with $Im\,\mu \neq 0$, $w(x) \in L^2(\mathbb{R}^2)$, and

$$\int_0^\infty g(\tau)\tau^{1+|m|}\,d\tau = 0,$$

then $w \equiv 0$.

Proof. Under the assumptions of the Lemma (1.6) implies

$$\psi(s) = -\frac{1}{2|m|}\,s^{|m|}\int_s^\infty g(\tau)\tau^{1-|m|}\,d\tau + \frac{1}{2|m|}\,s^{-|m|}\int_s^\infty g(\tau)\,\tau^{1+|m|}\,d\tau.$$

The same bootstrapping argument as above yields the existence of a constant $C_N$ for any $N \in \mathbb{Z}_+$, such that

$|g(s)| \leq C_N(1+s)^{-N}$, $\forall s \geq 0$.

Since $Im\,\mu \neq 0$, (1.11) implies

$|g(s)| \leq C(1+s)^{-2-\alpha}[s^{|m|}\int_s^\infty |g(\tau)|\,\tau^{|m|+2+\alpha}\,\tau^{-2|m|-1-\alpha}\,d\tau$

$+ s^{-|m|}\int_s^\infty |g(\tau)|\,\tau^{|m|+2+\alpha}\,\tau^{-1-\alpha}d\tau]. \tag{1.21}$

Define for any $s \geq 0$

$f(s) = \sup_{\tau \geq s}|g(\tau)|\,\tau^{|m|+2+\alpha}.$

The inequality (1.21) implies

$$f(s) \leq Cs^{-\alpha}f(s), \qquad \forall s > 0.$$

Thus $g(s) = 0$ for $s \geq C$, where the constant $C$ is sufficiently large. Gronwall's inequality, applied to (1.21) yields the statement of the lemma. QED.

We study the eigenfunction equation (1.11) with $Im\,\mu > 0$. From (1.5), (1.11)

$$-\frac{d^2}{ds^2}\psi - \frac{1}{s}\frac{d}{ds}\psi + \frac{m^2}{s^2}\psi + \frac{B(s)}{R(s)-\mu}\psi = 0. \qquad (1.22)$$

For the first observation about (1.22) lets introduce a new variable $t = \log s$, $t \in \mathbb{R}$.

Then (1.22) becomes

$$L\psi \equiv -\frac{d^2}{dt^2}\psi + m^2\psi + \frac{A(t)}{\Omega(t)-\mu}\psi = 0, \qquad (1.23)$$

where

$$A(t) = e^{2t}B(e^t), \qquad (1.24)$$

$$\Omega(t) = R(e^t). \qquad (1.25)$$

Therefore,

$$A(t) = \frac{d}{dt}G(e^t); \qquad (1.26)$$

$$\Omega(t) = \int_{-\infty}^{t} e^{-2(t-\tau)}G(e^\tau)d\tau; \qquad (1.27)$$

$$A(t) = \Omega''(t) + 2\Omega'(t), \quad (.)' = \frac{d}{dt}(.). \qquad (1.28)$$

From (1.2), (1.14)

$$A(t) = -\alpha e^{-\alpha t}, \quad t \geq \log M; \qquad (1.29)$$

$$\Omega(t) = C\alpha e^{-2t} + \frac{1}{2-\alpha}e^{-\alpha t}, \quad s \geq \log M; \qquad (1.30)$$

We will prove the existence of a nontrivial solution to (1.23), (1.28) with

$$Im\,\mu > 0 \qquad (1.31)$$

for a certain class of

potentials $\Omega(t)$ that we are going to describe now.

Class $\mathcal{C}$.

We denote by $\mathcal{C}$ the following class of functions

$$\Omega: \mathbb{R} \to (0, \infty)$$

$$(i)\ \Omega \in C^\infty(\mathbb{R}),\ \Omega'(t) < 0\ and$$

$$\int_{-\infty}^{t} e^{2\tau} A(\tau) d\tau < 0, \forall\, t \in \mathbb{R}; \tag{1.32}$$

$(ii)\ \Omega(t) = \Omega(-\infty) - c_0 e^{2t}, for\ some\ c_0 > 0, \qquad \forall t \leq \log M_1,\ where\ 0 < M_1 < M. \tag{1.33}$

This implies

$$G(s) = 2\Omega(-\infty) - 4c_0\, s^2, \forall\, s \in [0, M_1].$$

$(ii)\ A(t)\ has\ exactly\ 2\ zeroes, \qquad a < b,\ A'(a) > 0,\ A'(b) < 0.$

Remark 1. The inequality (1.32) holds if and only if

$$\int_0^{e^b} s^2 G'(s) ds < 0.$$

Remark 2. The inequality (1.32) implies

$$\Omega'(t) < 0\ \forall\, t \in \mathbb{R}.$$

Indeed, assuming (1.32), we have from (1.27)

$$\Omega'(t) = G(e^t) - 2 \int_{-\infty}^{t} e^{-2(t-\tau)} G(e^\tau) d\tau$$

$$= G(e^t) - e^{-2(t-\tau)}\, G(e^\tau)\,|_{\tau=-\infty}^{\tau=t} + \int_{-\infty}^{t} e^{-2(t-\tau)} G(e^\tau) d\tau$$

$$= e^{-2t} \int_{-\infty}^{t} e^{3\tau} G'(e^\tau) d\tau$$

$$= e^{-2t} \int_0^{e^t} s^2 G'(s) ds\ < 0. \tag{1.34}$$

For $m > 0$ we define the kernel

$$K_m(\xi, \eta) = \frac{1}{2m} e^{-m|\xi-\eta|}, \qquad where\ \xi, \eta \in \mathbb{R}. \tag{1.35}$$

This is the Green's function for the operator

$$-\frac{d^2}{dt^2} + m^2$$

on $\mathbb{R}$.

Indeed, for $f \in L^p(\mathbb{R}^2)$ with some $p \in [1, \infty]$, define

$$\psi(t) = \int K_m(\xi, \eta) f(\eta) \, d\eta$$

$$= e^{mt} \frac{1}{2m} \int_t^\infty e^{-m\eta} f(\eta) d\eta + e^{-mt} \frac{1}{2m} \int_{-\infty}^t e^{m\eta} f(\eta) d\eta.$$

Then,

$$-\psi''(t) + m^2 \psi(t) = -[\frac{1}{2} e^{mt} \int_t^\infty e^{-m\eta} f(\eta) d\eta - \frac{1}{2m} f(t) - \frac{1}{2} e^{-mt} \int_{-\infty}^t e^{m\eta} f(\eta) d\eta$$

$$+ \frac{1}{2m} f(t)]' + m^2 \psi(t)$$

$$= -\frac{m}{2} e^{mt} \int_t^\infty e^{-m\eta} f(\eta) d\eta + \frac{1}{2} f(t) - \frac{m}{2} e^{-mt} \frac{1}{2} \int_{-\infty}^t e^{m\eta} f(\eta) d\eta + \frac{1}{2} f(t) + m^2 \psi(t)$$

$$= f(t),$$

as claimed.

The eigenfunction equation (1.23) with $\psi \in L^p(\mathbb{R}^2)$, $p \in [1, \infty]$ can be written as follows:

$$\psi(t) + \int K_m(t, \xi) \frac{A(\xi)}{\Omega(\xi) - \mu} \psi(\xi) d\xi = 0.$$

Notice that $\psi \in W^{1,2}(\mathbb{R}, dt)$ is equivalent to $\partial_x \left( e^{im\theta} \psi(\log |x|) \right) \in L^2(\mathbb{R}^2)$, where $m \neq 0$ is an integer.

Indeed, the radial component of the gradient gives

$$\int_0^\infty s^{-2} |\psi'(\log s)|^2 s \, ds < \infty.$$

This integral equals

$$\int_{-\infty}^\infty |\psi'(t)|^2 dt.$$

The tangential component of the gradient gives

$$m^2 \int_0^\infty s^{-1} |\psi(\log s)|^2 \, ds < \infty.$$

After changing variables this yields

$$m^2 \int_{-\infty}^\infty |\psi(t)|^2 \, dt < \infty.$$

We define for every real $m > 0$ the set

$$\mathcal{U}_m = \{\mu \in \mathbb{C} | \, Im \, \mu > 0 \text{ and } \exists \, \psi = \psi_{m,\mu} \in L^2 \, (\mathbb{R}), \psi \neq 0 \, solving \, the \, equation \, (1.23)\}.$$

Remark. Any eigenfunction $\psi \in L^2 \, (\mathbb{R})$ of the problem (1.23) satisfies $\psi \in W^{1,2} \, (\mathbb{R})$.

Here is the essential technical result on linear instability of radially symmetric flows.

Theorem 1.1. (Theorem 11.1) For any $\alpha \in (0,2)$ there exists a function $\Omega \in \mathcal{C}$ and an integer $m \geq 2$, so that

$$\#\mathcal{U}_m = 1 \text{ and } \mathcal{U}_{m+l} = \emptyset \, for \, any \, positive \, integer \, l.$$

The proof of is given in §11 below.

## § 2. Asymptotic behavior of an eigenfunction

We need the following lemma.

Lemma 2.1. Let $v: \mathbb{R} \to \mathbb{R}$ be a measurable function so that $v \in L^1((-\infty, N))$ for every real $N$.

For $m > 0$ consider the equation

$$-\frac{d^2}{dt^2} y + m^2 y + v(t) y = 0.$$

There exists a unique solution $y(t)$, such that

$$y(t) \sim e^{mt}$$

as $t \to -\infty$. For this solution

$$y(t) = e^{mt}(1 + z(t)), \quad \forall t \in \mathbb{R},$$

where

$$|z(t)| \leq \exp\{\frac{1}{2m}\int_{-\infty}^{t}|v(s)|ds\} - 1, \forall t \in \mathbb{R}, \tag{2.1}$$

$$|z'(t)| \leq 2m \left(\exp\{\frac{1}{2m}\int_{-\infty}^{t}|v(s)|ds\} - 1\right), \forall t \in \mathbb{R}. \tag{2.2}$$

Proof. We have

$$y(t) - \frac{1}{m}\int_{-\infty}^{t}v(s)y(s)sh\left(m(t-s)\right)ds = 0 \; .$$

If $y(t) = e^{mt}(1 + z(t))$,

then

$$1 + z(t) - \frac{1}{2m}\int_{-\infty}^{t}v(s)\left(1 - e^{-2m(t-s)}\right)(1 + z(s))ds = 0.$$

Solving in $(1 + z(t))$ by iterations and comparing with the solution of the ODE

$$\frac{d}{dt}(1 + z_1(t)) = \frac{1}{2m}|v(t)|(1 + z_1(t)),$$

we arrive at the inequality

$$|z(t)| \leq \exp\{\frac{1}{2m}\int_{-\infty}^{t}|v(s)|ds\} - 1, \forall t \in \mathbb{R}.$$

This proves (2.1).

Also,

$$z'(t) - \int_{-\infty}^{t}v(s)e^{-2m(t-s)}(1 + z(s))ds = 0, \quad \forall t \in \mathbb{R}.$$

Using (2.1), we get

$$|z'(t)| \leq \int_{-\infty}^{t}|v(s)|\exp\{\frac{1}{2m}\int_{-\infty}^{s}|v(\xi)| \, d\xi\} \, ds$$

$$= 2m(exp\{\frac{1}{2m}\int_{-\infty}^{t}|v(s)|\,ds\}-1), \qquad \forall t \in \mathbb{R}.$$

This proves (2.2) and concludes the proof of Lemma 2.1. QED.

## § 3. Dispersion curves $(\mathcal{U}_m, m)$, $\mathcal{U}_m \neq \emptyset$.

We study the evolution of the set $\mathcal{U}_m$ as $m \searrow$ from $\infty$. The following 4 Propositions are of importance.

Proposition 3.1. The point $(m_0, 0)$ for $m_0 \geq 1$ is not a limit point of the set

$$\{(m, \mu) \in (0, \infty) \times \mathbb{C} \mid \mu \in \mathcal{U}_m\} \ .$$

Proof. Assume to the contrary, there is sequence $(m_j, \mu_j), j = 1,2,3,\ldots;\ m_j \to m_0 \geq 1,\ \mu_j \to 0$ as $j \to \infty$. Also $Im\ \mu_j > 0$ and $\mu_j \in \mathcal{U}_{m_j}$, $\forall j = 1,2,3,\ldots$. Using Lemma 2.1, we normalize $\psi_j(t) = \psi_{m_j,\mu_j}(t)$ so that

$$\psi_j(t) \sim e^{m_j t}$$

as $t \to -\infty$. Then

$$\psi_j(t) \sim C_j e^{-m_j t}$$

as $t \to \infty$ with an appropriate $C_j \in \mathbb{C}$.

We multiply both sides of (1.23) by $\overline{\psi_j(t)}$, integrate and take the imaginary part. This implies (here we use $Im\ \mu_j > 0$)

$$\int_{-\infty}^{\infty} \frac{\Omega'' + 2\Omega'}{|\Omega - Re\ \mu_j|^2 + |Im\ \mu_j|^2}|\psi_j|^2\,dt = 0. \tag{3.1}$$

We split the integral in the left side of (3.1) as follows:

$$\int_{-\infty}^{a} + \int_{a}^{b} + \int_{b}^{\infty} = 0. \tag{3.2}$$

We have for every $N \in \mathbb{R}$, $\psi_j \to \psi_0$ uniformly on $(-\infty, N)$ as $j \to \infty$, where $\psi_0$ solves the equation

$$-\frac{d^2}{dt^2}\psi_0 + m_0^2 \psi_0 + \Omega^{-1}(\Omega'' + 2\Omega')\psi_0 = 0, \ t \in \mathbb{R}, \tag{3.3}$$

$$\psi_0(t) \sim e^{m_0 t},$$

as $t \to -\infty$. From (3.1), (3.2) using Lebesgue dominated convergence theorem we pass to the limit as $j \to \infty$ to obtain for any fixed $N \geq b$

$$-\int_{-\infty}^{a} \Omega^{-2}(\Omega'' + 2\Omega')|\psi_0|^2 \, dt - \int_{b}^{N} \Omega^{-2}(\Omega'' + 2\Omega')|\psi_0|^2 \, dt$$

$$\leq \int_{a}^{b} \Omega^{-2}(\Omega'' + 2\Omega')|\psi_0|^2 \, dt \ .$$

Passing to the limit as $N \to \infty$ (note the integrands in the left side are nonpositive) we get

$$\Omega^{-2}(\Omega'' + 2\Omega')|\psi_0|^2 \in L^1(\mathbb{R}). \tag{3.4}$$

By (1.29), (1.30) $\Omega^{-1}(\Omega'' + 2\Omega') \to -\alpha(2-\alpha)$ as $t \to \infty$. Therefore, the only nontrivial solution that does not violate (3.4) is the one satisfying

$$\psi_0(t) \sim C_1 e^{-\sqrt{m_0^2 - \alpha(2-\alpha)}\, t} \ \text{as } t \to \infty,$$

where $C_1 \in \mathbb{R}$, $C_1 \neq 0$.

See, e.g., [BS], p.61.

The only exceptional case is $m_0 = 1, \alpha = 1$, which is easily ruled out by (3.4).

Using (3.4) again we get

$$m_0^2 > \frac{3}{4}\alpha\left(\frac{8}{3} - \alpha\right) \ . \tag{3.5}$$

If (3.5) fails this is a contradiction, therefore we assume (3.5). From (3.3)

$$\int |\frac{d}{dt}\psi_0 - \frac{\Omega'}{\Omega}\psi_0|^2 + m_0^2|\psi_0|^2 + 2\frac{\Omega'}{\Omega}|\psi_0|^2 \, dt = 0 \ , \tag{3.6}$$

where the integration by parts is legitimate.

Let $\psi_0 = \Omega\eta$, $\eta(t) \sim \Omega(-\infty)^{-1} e^{m_0 t}$ as $t \to -\infty$; $\eta(t) \sim C_1(2-\alpha)e^{-\left(\sqrt{m_0^2 - \alpha(2-\alpha)} - \alpha\right)t}$, as $t \to \infty$.

We get from (3.6) after integration by parts

$$\int \Omega^2 \left|\frac{d\eta}{dt}\right|^2 + m_0^2 \Omega^2 |\eta|^2 + 2\Omega\Omega' |\eta|^2 dt = 0,$$

therefore,

$$\int \Omega^2 (\left|\frac{d\eta}{dt}\right|^2 + m_0^2 |\eta|^2 - 2\eta\frac{d\eta}{dt}) dt = 0.$$

If $m_0 > 1$, this implies $\eta \equiv 0$, a contradiction. If $m_0 = 1$, this implies $\eta = Ce^t$, in contradiction with (3.4). This concludes the proof of Proposition 3.1. QED.

Proposition 3.2. For any $m_0 > 0$  $(m_0, \Omega(-\infty))$ is not a limit point of the set

$$\{(m, \mu) \in (0, \infty) \times \mathbb{C} \mid \mu \in \mathcal{U}_m\} .$$

Proof. Let to the contrary

$$\psi_j(t) = \psi_{m_j, \mu_j}(t), \quad j = 1,2,3 \ldots$$

satisfy (1.23) with $m = m_j$, $\mu = \mu_j$; $m_j \to m_0 > 0$, $\mu_j \to \Omega(-\infty)$ as $j \to \infty$; $Im\, \mu_j > 0$ for $j = 1,2,3 \ldots$ .

We assume normalization

$$\psi_j(t) \sim e^{-m_j t}, \quad t \to \infty. \tag{3.7}$$

We have, using $Im\, \mu_j > 0$

$$\int_{-\infty}^{\infty} \frac{\Omega'' + 2\Omega'}{|\Omega - Re\, \mu_j|^2 + |Im\, \mu_j|^2} |\psi_j|^2 \, dt = 0. \tag{3.8}$$

On every interval $[N, \infty)$, $N \in \mathbb{R}$, $\psi_j \to \psi_0$ uniformly as $j \to \infty$, where $\psi_0$ satisfies the Sturm-Liouville equation on the line

$$-\frac{d^2}{dt^2}\psi_0 + m_0^2 \psi_0 + (\Omega - \Omega(-\infty))^{-1}(\Omega'' + 2\Omega')\psi_0 = 0, \quad t \in \mathbb{R},$$

$$\psi_0(t) \sim e^{-m_0 t}, \quad t \to \infty . \tag{3,9}$$

Obviously, $\psi_0(t)$ is real valued.

From (3.8), passing to the limit as $j \to \infty$,

$$-\int_b^\infty (\Omega - \Omega(-\infty))^{-2} (\Omega'' + 2\Omega')|\psi_0|^2 \, dt - \int_N^a (\Omega - \Omega(-\infty))^{-2} (\Omega'' + 2\Omega')|\psi_0|^2 \, dt$$

$$\leq \int_a^b (\Omega - \Omega(-\infty))^{-2} (\Omega'' + 2\Omega')|\psi_0|^2 \, dt \, , \quad \forall \, N \leq a \, . \tag{3.10}$$

Passing to the limit as $N \to -\infty$ yields

$$(\Omega - \Omega(-\infty))^{-2}(\Omega'' + 2\Omega')|\psi_0|^2 \in L^1(\mathbb{R}) \, . \tag{3.11}$$

From (1.33)

$$(\Omega - \Omega(-\infty))^{-1}(\Omega'' + 2\Omega') = 8 \, , \quad \forall t \leq \log M_1.$$

Therefore,

$$\psi_0(t) = C_1 e^{-t\sqrt{m_0^2+8}} + C_2 e^{t\sqrt{m_0^2+8}} \, , \quad \forall t \leq \log M_1.$$

To satisfy (3.11) this solution must have $C_1 = 0$. From (3.9) integrating by part yields

$$\int |\frac{d}{dt}\psi_0|^2 + m_0^2|\psi_0|^2 + (\Omega - \Omega(-\infty))^{-2}\Omega'^2|\psi_0|^2 - 2(\Omega - \Omega(-\infty))^{-1}\Omega' \, \psi_0' \psi_0$$

$$+ 2(\Omega - \Omega(-\infty))^{-1}\Omega'|\psi_0|^2 \, dt = 0 \, .$$

Therefore,

$$\int |\frac{d}{dt}\psi_0 - (\Omega - \Omega(-\infty))^{-1}\Omega' \, \psi_0|^2 + m_0^2|\psi_0|^2 + 2(\Omega - \Omega(-\infty))^{-1}\Omega'|\psi_0|^2 \, dt = 0.$$

But $(\Omega - \Omega(-\infty))^{-1}\Omega' > 0$, $\forall t \in \mathbb{R}$. Therefore, $\psi_0 \equiv 0$ in contradiction with construction. This proves Proposition 3.2. QED.

Proposition 3.3. Let $(m_0, \mu_0)$, where $m_0 \in [1, \infty)$, $\mu_0 \in (0, \Omega(-\infty))$ be a limit point of the set

$$\{(m, \mu) \in (0, \infty) \times \mathbb{C} \mid \mu \in \mathcal{U}_m\} \, . \tag{3.12}$$

Then either $\mu_0 = \Omega(a)$ or $\mu_0 = \Omega(b)$. Also, spectral problem

$$-\frac{d^2}{dt^2}\psi_0 + m_0{}^2\psi_0 + (\Omega(t) - \mu_0)^{-1}A(t)\psi_0 = 0, \ t \in \mathbb{R}$$

has a solution $\psi_0 \in L^2(\mathbb{R})$, $\psi_0 \not\equiv 0$. In fact, $\psi_0$ decays exponentially at $\pm\infty$.

Proof. Let to the contrary, a point $(m_0, \mu_0)$, where $m_0 \in (0, \infty)$, $\mu_0 \in \big(0, \Omega(-\infty)\big)\backslash(\{\Omega(a)\} \cup \{\Omega(b)\})$ be the limit point of the set (3.12). Let $\mu_0 = \Omega(c_0)$, $c_0 \neq a$, $c_0 \neq b$. Then there is a sequence $(m_j, \mu_j) \to (m_0, \mu_0), j = 1,2,3, \ldots$ and for every $j$ an eigenfunction $\psi_j \in \psi_0 \in L^2(\mathbb{R})$ such that

$$-\frac{d^2}{dt^2}\psi_j + m_j{}^2\psi_j + \big(\Omega(t) - \mu_j\big)^{-1}A(t)\psi_j = 0, \ , \ \text{Im}\ \mu_j > 0, \ t \in \mathbb{R}. \tag{3.13}$$

We use the following normalization of $\psi_j$:

$$\int |\frac{d}{dt}\psi_j|^2 + m_j{}^2|\psi_j|^2 \, dt = 1 \ ; \tag{3.14}$$

$$\psi_j(t) \sim \rho_j^{\pm} e^{\mp m_j t}, \quad t \to \pm\infty\ ; \tag{3.15}$$

$$\rho_j^+ > 0, \quad \rho_j^- \in \mathbb{C}\backslash\{0\}, \quad j = 1,2,3, \ldots \ . \tag{3.16}$$

By (3.15). (3.16) and Lemma 2.1 both sequences $\{\rho_j^+\}$ and $\{\rho_j^-\}$ are bounded (otherwise, (3.14) is violated). We can and will by possibly selecting a subsequence assume that

$$\psi_j \to \psi_0 \in L^2(\mathbb{R}), \ \psi_j' \to \psi_0' \in L^2(\mathbb{R}) \ weakly \ in \ L^2(\mathbb{R}). \tag{3.17}$$

Therefore, $\psi_0 \in W^{1,2}(\mathbb{R})$.

From Lemma 2.1 and using the asymptotic behavior of $\Omega(t)$ as $t \to \pm\infty$,

$$\forall \varepsilon > 0 \ \exists N \in \mathbb{R}, \ \int_{|t| \geq N}|\psi_j|^2 \, dt < \varepsilon\ ; \ j = 1,2,3, \ldots \ . \tag{3.18}$$

By Sobolev embedding $W^{1,2}(\mathbb{R}) \subset C^{\frac{1}{2}}(\mathbb{R})$

$$\|\psi_j\|_{L^\infty(\mathbb{R})} \leq C < \infty, \ \|\psi_j\|_{C^{\frac{1}{2}}(\mathbb{R})} \leq C < \infty\ ; j = 1,2,3, \ldots \ . \tag{3.19}$$

From (3.14), (3.18), (3.19) by possibly selecting a subsequence $\{\psi_j\}$, $j = 1,2,3, \ldots$ we may assume that

$$\psi_j \to \psi_0 \text{ strongly in } L^2(\mathbb{R}). \tag{3.20}$$

From (3.13), (3.14)

$$\int (\Omega(t) - \mu_j)^{-1} A(t) |\psi_j(t)|^2 \, dt = -1. \tag{3.21}$$

Therefore, for any $\delta > 0$

$$Im \int_{\mathbb{R}\setminus(c_0-\delta,c_0+\delta)} (\Omega(t) - \mu_j)^{-1} A(t) |\psi_j(t)|^2 \, dt$$

$$+Im \int_{(c_0-\delta,c_0+\delta)} (\Omega(t) - \mu_j)^{-1} A(t)(|\psi_j(t)|^2 - |\psi_j(c_j)|^2 + |\psi_j(c_j)|^2) \, dt = 0. \tag{3.22}$$

Here $\mu_j = c_j + id_j$, $d_j > 0$, $j = 1,2,3,\ldots$. From (3.19), (3.20) $\psi_j \to \psi_0$ uniformly on $\mathbb{R}$, therefore, since $\mu_j \to \mu_0 = \Omega(c_0) \in \mathbb{R}$, the first term in the left side of (3.22) converges to 0 as $j \to \infty$. Since $\|\psi_j\|_{C^{\frac{1}{2}}(\mathbb{R})} \leq C < \infty$, where the constant $C$ is uniform with respect to $j$, we have

$$Im \int_{(c_0-\delta,c_0+\delta)} (\Omega(t) - \mu_j)^{-1} A(t)(|\psi_j(t)|^2 - |\psi_j(c_j)|^2) \, dt = O\left(\delta^{\frac{1}{2}}\right) \tag{3.23}$$

with the constant in $O$ uniform with respect to $j$. Passing to the limit as $j \to \infty$ in (3.22) and using (3.23) and the argument above, we get

$$\pi |\psi_0(c_0)|^2 A(c_0) |\Omega'(c_0)|^{-1} = 0.$$

Since $c_0 \neq a, c_0 \neq b$, we have

$$\psi_0(c_0) = 0. \tag{3.24}$$

Therefore, in the sense of $S'(\mathbb{R})$

$$(\Omega(t) - \mu_j)^{-1} \psi_j \to (\Omega(t) - \mu_0)^{-1} \psi_0.$$

Passing to the limit as $j \to \infty$ in the equation (3.13) we get

$$-\frac{d^2}{dt^2}\psi_0 + m_0^2 \psi_0 + (\Omega(t) - \mu_0)^{-1} A(t)\psi_0 = 0, \quad t \in \mathbb{R} \tag{3.25}$$

in $S'(\mathbb{R})$. From (3.19)

$$(\Omega(t) - \mu_0)^{-1}\psi_0 \in L^{2-\theta}(\mathbb{R}) \ , \forall \theta \in (0,1).$$

It follows from the equation (3.25) and from boundedness of the pseudo-differential operator with symbol $(|\xi|^2 + m_0^2)^{-1}|\xi|^2$ in $L^{2-\theta}(\mathbb{R})$, that $\frac{d^2}{dt^2}\psi_0 \in L^{2-\theta}(\mathbb{R})$ for any $\theta \in (0,1)$.

By Sobolev embedding theorem this yields $\frac{d}{dt}\psi_0 \in C^\sigma(\mathbb{R})$ for any $\sigma \in \left(0, \frac{1}{2}\right)$. Obviously, $\psi_0$ is real valued on $(c_0, \infty)$,

Integrating by parts and using $\psi_0(t) = O(|t - c_0|)$ as $t \to c_0$ to justify it, we arrive at

$$\int_{c_0}^\infty |\frac{d}{dt}\psi_0|^2 + m_0^2|\psi_0|^2 + (\Omega - \mu_0)^{-2}\Omega'^2|\psi_0|^2 - 2(\Omega - \mu_0)^{-1}\Omega'\,\psi_0\,\psi_0'$$

$$+ 2(\Omega - \mu_0)^{-1}\Omega'|\psi_0|^2\,dt = 0 \ .$$

Therefore,

$$\int_{c_0}^\infty |\frac{d}{dt}\psi_0 - (\Omega - \mu_0)^{-1}\Omega'\,\psi_0|^2 + m_0^2|\psi_0|^2 + 2(\Omega - \mu_0)^{-1}\Omega'|\psi_0|^2\,dt = 0.$$

But $(\Omega(t) - \mu_0)^{-1}\Omega'(t) > 0$ for $t \in (c_0, \infty)$. Therefore, $\psi_0 \equiv 0$ on $(c_0, \infty)$, and consequently

$$\psi_0'(c_0) = 0.$$

Thus $\psi_0(t) = O(|t - c_0|^{1+\sigma})$ as $t \to c_0$.

Likewise, legitimately integrating by parts on $\mathbb{R}$, we get

$$\int_{-\infty}^\infty |\frac{d}{dt}\psi_0 - (\Omega - \mu_0)^{-1}\Omega'\,\psi_0|^2 + m_0^2|\psi_0|^2 + 2(\Omega - \mu_0)^{-1}\Omega'|\psi_0|^2\,dt = 0.$$

Substituting $\psi_0(t) = (\Omega - \mu_0)\eta(t)$, we arrive at

$$\int (\Omega - \mu_0)^2 \left(\left|\frac{d\eta}{dt}\right|^2 + m_0^2|\eta|^2 - 2\mathrm{Re}\,(\bar{\eta}\frac{d\eta}{dt})\right)dt = 0.$$

If $m_0 > 1$ this implies $\eta(t) \equiv 0, \forall t \in \mathbb{R}$. If $m_0 = 1$ this implies $\eta(t) = Ce^t, \forall t \in \mathbb{R}$, where $C \in \mathbb{R}$. But $\eta(c_0) = 0$. Therefore, $C = 0$. In both cases, $\psi_0(t) \equiv 0$ on $\mathbb{R}$.

Therefore,

$$\psi_j \to 0 \text{ uniformly on } \mathbb{R} \text{ as } j \to \infty. \tag{3.26}$$

From this using Lemma 2.1 we conclude

$$\rho_j^\pm \to 0 \text{ as } j \to \infty. \tag{3.27}$$

Multiplying both sides of (3.13) by $\overline{(\psi_j(t) - \psi_j(c_j))}$ where $c_j = \operatorname{Re} \mu_j$ and integrating by parts we get

$$\int |\tfrac{d}{dt}\psi_j|^2 + m_j^2 \psi_j \overline{(\psi_j(t) - \psi_j(c_j))} + (\Omega(t) - \mu_j)^{-1} A(t)\psi_j(t) \overline{(\psi_j(t) - \psi_j(c_j))} \, dt = 0. \tag{3.28}$$

The integral of the second term goes to zero as $j \to \infty$ since $\psi_j(c_j) \to 0$ by (3.26), (3.27) and Lemma 2.1. The integral of the last term equals

$$\int_{\mathbb{R}\setminus(c_0-\delta,c_0+\delta)} (\Omega(t) - \mu_j)^{-1} A(t)\psi_j \overline{(\psi_j(t) - \psi_j(c_j))} \, dt$$

$$+ \int_{(c_0-\delta,c_0+\delta)} (\Omega(t) - \mu_j)^{-1} A(t)\psi_j \overline{(\psi_j(t) - \psi_j(c_j))} \, dt.$$

The first integral converges to zero for a fixed $\delta > 0$ as $j \to \infty$. The second integral is $O(\delta)$ with a constant uniform with respect to $j$. Indeed,

$$\left| (\Omega(t) - \mu_j)^{-1} \overline{(\psi_j(t) - \psi_j(c_j))} \right| \leq C \left( - \inf_{(c_0-\delta,c_0+\delta)} \Omega' \right)^{-1}.$$

For any fixed $\delta > 0$ passing to the limit in (3.28) we get

$$\limsup_{j\to\infty} \int |\tfrac{d}{dt}\psi_j|^2 \, dt \leq C\delta.$$

Therefore,

$$\lim_{j\to\infty} \int |\frac{d}{dt}\psi_j|^2 \, dt = 0.$$

But from (3.14), (3.20), with $\psi_0 = 0$ this limit is 1. Contradiction.

Therefore, either $\mu_0 = \Omega(a)$ or $\mu_0 = \Omega(b)$. Let $d = a$ or $d = b$ so that $\mu_0 = \Omega(d)$. Assume the problem

$$-\frac{d^2}{dt^2}\psi_0 + m_0{}^2\psi_0 + (\Omega(t) - \mu_0)^{-1}A(t)\psi_0 = 0, \ t \in \mathbb{R} \tag{3.29}$$

*has the only solution* $\psi_0 \in L^2(\mathbb{R})$, $\psi_0 \equiv 0$.

We have (3.13) satisfied with the normalization (3.14), (3.15), (3.16). Therefore, (3.17) --(3.20) hold true.

Passing to the limit as $j \to \infty$ in (3.13), we get

$$-\frac{d^2}{dt^2}\psi_0 + m_0{}^2\psi_0 + (\Omega(t) - \mu_0)^{-1}A(t)\psi_0 = 0, \ t \in \mathbb{R} \ . \tag{3.30}$$

To handle the last term in (3.13) we use the following representation

$$\left(\Omega(t) - \mu_j\right)^{-1} A(t)\psi_j(t) = \left(\Omega(t) - \mu_j\right)^{-1}(A(t) - A(c_j))\psi_j(t)$$

$$+ A(c_j)\left(\Omega(t) - \mu_j\right)^{-1}\left(\psi_j(t) - \psi_j(c_j)\right) + A(c_j)\psi_j(c_j)\left(\Omega(t) - \mu_j\right)^{-1}.$$

In $S'(\mathbb{R})$ the first term in the right side converges to $(\Omega(t) - \mu_0)^{-1}A(t)\psi_0(t)$ and the second and the third term both converge to 0 (since $A(c_j) \to 0$) as $j \to \infty$. Also,

$$\left(\Omega(t) - \mu_j\right)^{-1} \to (\Omega(t) - \mu_0 - i0)^{-1}$$

in $S'(\mathbb{R})$ as $j \to \infty$.

From (3.29), (3.30) it follows that $\psi_0 \equiv 0$. Therefore, from Lemma 2.1 (3.27) holds true. Multiplying both sides in (3.13) by $\overline{\psi_j(t)} - \overline{\psi_j(c_j)}$ and (legitimately) integrating by parts, we get

$$\int |\frac{d}{dt}\psi_j|^2 + m_j{}^2\psi_j\left(\overline{\psi_j(t)} - \overline{\psi_j(c_j)}\right) + (\Omega(t) - \mu_j)^{-1}A(t)\psi_j(t)\left(\overline{\psi_j(t)} - \overline{\psi_j(c_j)}\right) dt = 0. \tag{3.31}$$

The integral of the second term goes to zero as $j \to \infty$ since $\psi_j(c_j) \to 0$ by (3.26), (3.27) and Lemma 2.1. The integral of the last term equals

$$\int (\Omega(t) - \mu_j)^{-1} \left(A(t) - A(c_j)\right) \psi_j(t) \left(\overline{\psi_j(t)} - \overline{\psi_j(c_j)}\right) dt + o(1) \text{ as } j \to \infty.$$

We have

$$\left|(\Omega(t) - \mu_j)^{-1} \left(A(t) - A(c_j)\right)\right| \leq C < \infty, \forall t \in \mathbb{R}$$

with a constant $C$ independent of $j = 1,2,3,\ldots$ . Also $\|\psi_j\|_{L^2(\mathbb{R})} \to 0$; $\|\psi_j\|_{L^1(\mathbb{R})} \to 0$; $\|\psi_j\|_{C(\mathbb{R})} \to 0$ as $j \to \infty$ because of the uniform convergence $\psi_j \to \psi_0$ on $\mathbb{R}$ and because of (3.27). From (3.31) we get

$$\int \left|\tfrac{d}{dt}\psi_j\right|^2 + m_j^2|\psi_j|^2 \, dt \to 0$$

as $j \to \infty$. This contradicts our normalization (3.14). This concludes the proof of Proposition 3.3. QED.

Proposition 3.4. Let $m_0 \in [1, \infty)$, $\mu_0 \in \mathbb{R} \setminus [0, \Omega(-\infty)]$. Then $(m_0, \mu_0)$ is not a limit point of the set

$$\{(m, \mu) \in (0, \infty) \times \mathbb{C} \mid \mu \in \mathcal{U}_m\}.$$

Proof. Assume to the contrary, there is sequence $(m_j, \mu_j), j = 1,2,3,\ldots$; $m_j \to m_0 \geq 1$, $\mu_j \to \mu_0$ as $j \to \infty$. Also $\operatorname{Im} \mu_j > 0$ and $\mu_j \in \mathcal{U}_{m_j}$, $\forall j = 1,2,3,\ldots$ . Using Lemma 2.1, we normalize $\psi_j(t) = \psi_{m_j, \mu_j}(t)$ so that

$$\psi_j(t) \sim e^{m_j t}$$

as $t \to -\infty$. Then

$$\psi_j(t) \sim C_j e^{-m_j t}$$

as $t \to \infty$ with an appropriate $C_j \in \mathbb{C}$.

The eigenfunction $\psi_j(t)$ satisfies the following equation:

$$-\tfrac{d^2}{dt^2}\psi_j + m_j^2 \psi_j + (\Omega(t) - \mu_j)^{-1} A(t)\psi_j = 0, \quad \operatorname{Im} \mu_j > 0, \quad t \in \mathbb{R}. \tag{3.32}$$

Multiplying both sides in (3.32) by $\overline{\psi_j}(t)$ and separating the imaginary part we get (here we also use $Im\ \mu_j > 0$):

$$\int_{-\infty}^{\infty} \frac{\Omega'' + 2\Omega'}{|\Omega - Re\ \mu_j|^2 + |Im\ \mu_j|^2} |\psi_j|^2\ dt = 0. \tag{3.33}$$

Therefore,

$$-\int_{-\infty}^{a} (|\Omega(t) - Re\ \mu_j|^2 + |Im\ \mu_j|^2)^{-1} (\Omega''(t) + 2\Omega'(t)) |\psi_j(t)|^2\ dt$$

$$-\int_{b}^{\infty} (|\Omega(t) - Re\ \mu_j|^2 + |Im\ \mu_j|^2)^{-1} (\Omega''(t) + 2\Omega'(t)) |\psi_j(t)|^2\ dt$$

$$\leq \int_a^b (|\Omega(t) - Re\ \mu_j|^2 + |Im\ \mu_j|^2)^{-1} (\Omega''(t) + 2\Omega'(t)) |\psi_j(t)|^2\ dt\ . \tag{3.34}$$

Both integrals in the left side of (3.34) are positive.

From Lemma 2.1 and from the inequality (3.34) we conclude that

$$|C_j| \leq C < \infty, \forall j = 1,2,3, \dots\ . \tag{3.35}$$

Using again Lemma 2.1, (3.35) after possibly extracting a subsequence, we can pass to the limit in (3.32) as $j \to \infty$ to get

$$\psi_j \to \psi_0\ in\ S'(\mathbb{R}), \qquad \psi_0 \in L^2(\mathbb{R});$$

$$-\frac{d^2}{dt^2}\psi_0 + m_0^2 \psi_0 + (\Omega(t) - \mu_0)^{-1}(\Omega''(t) + 2\Omega'(t))\psi_0 = 0,\ ,\quad t \in \mathbb{R}; \tag{3.36}$$

$$\psi_0(t) \sim e^{m_0 t}\ as\ t \to -\infty;\ \psi_0(t) \sim \rho e^{-m_0 t}\ as\ t \to \infty.$$

Here $\rho \in \mathbb{R}, \rho \neq 0$. The eigenfunction $\psi_0$ is real-valued. Multiplying both sides in (3.36) by $\psi_0(t)$ and legitimately integrating by parts yields

$$\int |\frac{d}{dt}\psi_0|^2 + m_0^2 |\psi_0|^2 + (\Omega - \mu_0)^{-2} \Omega'^2 |\psi_0|^2 - 2(\Omega - \mu_0)^{-1} \Omega'\ \psi_0'\psi_0$$

$$+ 2(\Omega - \mu_0)^{-1} \Omega' |\psi_0|^2\ dt = 0\ .$$

Therefore,

$$\int |\frac{d}{dt}\psi_0 - (\Omega - \mu_0)^{-1} \Omega'\ \psi_0|^2 + m_0^2 |\psi_0|^2 + 2(\Omega - \mu_0)^{-1} \Omega' |\psi_0|^2\ dt = 0. \tag{3.37}$$

Using the substitution

$$\psi_0(t) = (\Omega(t) - \mu_0)\eta_0(t),$$

where $\eta_0 \in C^\infty(\mathbb{R})$ decays exponentially at $\pm\infty$, we get

$$\int |\Omega - \mu_0|^2 |\tfrac{d}{dt}\eta_0|^2 + m_0^2|\Omega - \mu_0|^2|\eta_0|^2 + 2\,(\Omega - \mu_0)'(\Omega - \mu_0)|\eta_0|^2\, dt = 0. \tag{3.38}$$

Integrating by parts in the right side of (3.38) yields

$$\int |\Omega - \mu_0|^2 (|\tfrac{d}{dt}\eta_0|^2 + m_0^2|\eta_0|^2 - 2\eta_0'\eta_0)\, dt = 0.$$

For $m_0 > 1$ this identity implies $\psi_0 \equiv 0$ in contradiction with construction. If $m_0 = 1$ this implies $\psi_0(t) = (\Omega(-\infty) - \mu_0)^{-1}(\Omega(t) - \mu_0)e^t$, $\forall t \in \mathbb{R}$. This contradicts the construction (the asymptotic behavior of $\psi_0(t)$ as $t \to \infty$.) This concludes the proof of Proposition 3.4. QED.

## § 4. The limiting equation

We next discuss the case $\mu \to \Omega(d)$, where $d = a$ or $d = b$. In this case, under certain conditions, a nontrivial branch of eigenvalues with $Im\,\mu > 0$ does appear. First, we study the limiting equation

$$-\frac{d^2}{dt^2}\psi_0 + m_0^2\psi_0 + (\Omega(t) - \Omega(d))^{-1}A(t)\psi_0 = 0, \ \ t \in \mathbb{R}, \tag{4.1}$$

Where $d = a$ or $d = b$. In this case the potential

$$v(t) = (\Omega(t) - \Omega(d))^{-1}A(t) \in S(\mathbb{R}).$$

Proposition 4.1. For a fixed $d = a$ or $d = b$, there exists at most one $m_0 \geq 1$ so that (4.1) has a nontrivial solution $\psi_0 \in L^2(\mathbb{R})$.

Proof. Consider

$$L \equiv -\frac{d^2}{dt^2} + v(t).$$

Let, to the contrary, there exist $m_1 > m_2 \geq 1$ so that (4.1) holds with $m_0 = m_{1,2}$, $\psi_0 = \psi_{1,2}$. The eigenfunctions $\psi_{1,2}$ decay exponentially at $\pm\infty$ (see Lemma 2.1). We may and will assume $\psi_{1,2}$ are real valued. There is a nontrivial linear combination $\psi = C_1\psi_1 + C_2\psi_2$ so that $\psi(d) = 0$. From (4.1)

$$\int |\frac{d}{dt}\psi|^2 + (\Omega(t) - \mu_0)^{-1}A(t)|\psi|^2 dt$$

$$= -C_1^2\, m_1^2 \int |\psi_1|^2\, dt - C_2^2\, m_2^2 \int |\psi_2|^2\, dt.$$

Therefore, integrating by parts and using that $\psi(d) = 0$, $\psi \in S(\mathbb{R})$

$$\int |\frac{d}{dt}\psi|^2 + (\Omega - \mu_0)^{-2}(\Omega' + 2\Omega)\Omega'|\psi|^2 - 2(\Omega - \mu_0)^{-1}(\Omega' + 2\Omega)\psi\psi' dt$$

$$= -C_1^2\, m_1^2 \int |\psi_1|^2\, dt - C_2^2\, m_2^2 \int |\psi_2|^2\, dt.$$

This implies

$$\int |\frac{d}{dt}\psi - (\Omega - \mu_0)^{-1}\Omega'\psi|^2 - 2\Omega \frac{d}{dt}[(\Omega - \mu_0)^{-1}]|\psi|^2 - 4\Omega(\Omega - \mu_0)^{-1}\psi\psi' dt$$

$$= -C_1^2\, m_1^2 \int |\psi_1|^2\, dt - C_2^2\, m_2^2 \int |\psi_2|^2\, dt.$$

Integrating by parts again and using $\psi(d) = 0$ yields

$$\int |\frac{d}{dt}\psi - (\Omega - \mu_0)^{-1}\Omega'\psi|^2 + 2(\Omega - \mu_0)^{-1}\Omega'\,|\psi|^2\, dt$$

$$= -C_1^2\, m_1^2 \int |\psi_1|^2\, dt - C_2^2\, m_2^2 \int |\psi_2|^2\, dt.$$

Let $\psi(t) = (\Omega(t) - \mu_0)\eta(t)$, where $\eta \in S(\mathbb{R})$. Then

$$\int |(\Omega(t) - \mu_0)\frac{d}{dt}\eta|^2 + 2(\Omega - \mu_0)\,\Omega'\,|\eta|^2\, dt$$

$$= -C_1^2\, m_1^2 \int |\psi_1|^2\, dt - C_2^2\, m_2^2 \int |\psi_2|^2\, dt.$$

This implies

$$\int |(\Omega(t) - \mu_0)|^2 |\frac{d}{dt}\eta - \eta|^2\, dt$$

$$= -C_1^2\, (m_1^2 - 1) \int |\psi_1|^2\, dt - C_2^2\, (m_2^2 - 1) \int |\psi_2|^2\, dt.$$

Since the right side is nonpositive, while the left side is strictly positive, this is a contradiction.

This completes the proof of Proposition 4.1. QED.

Remark. A double eigenvalue $m_0^2 \geq 1$ in (4.1) not possible by the uniqueness statement in Lemma 2.1.

Also, for an eigenfunction $\psi \in L^2(\mathbb{R})$ with an eigenvalue $m_0^2$ such that $\psi(d) = 0$ integration by parts yields

$$\int_d^\infty |\frac{d}{dt}\psi - (\Omega - \Omega(d))^{-1}\Omega'\psi|^2 + m_0^2|\psi|^2 + 2(\Omega - \Omega(d))^{-1}\Omega'|\psi|^2 \, dt = 0.$$

This leads to contradiction since $(\Omega(t) - \Omega(d))^{-1}\Omega'(t) > 0 \,, \forall t \in (d, \infty)$.

Indeed. $\psi(d) = 0, \psi'(d) = 0$ imply $\psi \equiv 0$. QED.

Proposition 4.2. For $d = a$ or for $d = b$, there exists a $m_0 > 1$ so that (4.1) has a non-trivial solution $\psi_0 \in L^2(\mathbb{R})$.

Proof. We consider

$$L \equiv -\frac{d^2}{dt^2} + v(t), \quad v(t) = (\Omega(t) - \Omega(d))^{-1}A(t), \quad d = a \text{ or } d = b,$$

as (unbounded) self-adjoint operator in $L^2(\mathbb{R})$. We compute the bottom of the spectrum of $L$ using the variational formula

$$-m_0^2 = \inf_{u \in W^{1,2}(\mathbb{R})\setminus\{0\}} \frac{1}{\|u\|_{L^2(\mathbb{R})}^2} \int \left|\frac{d}{dt}u\right|^2 + v(t)|u|^2 dt. \tag{4.2}$$

For $d = b$ we use

$$u(t) = \begin{cases} (\Omega(t) - \Omega(b))e^t, & t \leq b \\ 0, & t > b \end{cases}$$

in (4.2). With this choice,

$$\int \left|\frac{d}{dt}u\right|^2 + v(t)|u|^2 dt$$

$$= \int_{-\infty}^b -u''u + v(t)|u|^2 dt$$

$$= \int_{-\infty}^{b} e^{2t}(-\Omega'' - 2\Omega' - (\Omega - \Omega(b)))(\Omega - \Omega(b)) + e^{2t}(\Omega - \Omega(b))^{-1}(\Omega" + 2\Omega')(\Omega - \Omega(b))^2 dt$$

$$= -\int_{-\infty}^{b} e^{2t}(\Omega - \Omega(b))^2 dt = -\int |u|^2 dt.$$

Therefore, the right side in (4.2) is $< -1$.

For $d = a$ we use

$$u(t) = \begin{cases} (\Omega(t) - \Omega(a))e^t, & t \le a \\ 0, & t > a \end{cases}$$

in (4.2). With this choice,

$$\int \left|\frac{d}{dt}u\right|^2 + v(t)|u|^2 dt$$

$$= \int_{-\infty}^{a} -u''u + v(t)|u|^2 dt$$

$$= \int_{-\infty}^{a} e^{2t}(-\Omega'' - 2\Omega' - (\Omega - \Omega(a))(\Omega - \Omega(a)) + e^{2t}(\Omega - \Omega(a))^{-1}(\Omega" + 2\Omega')(\Omega - \Omega(a))^2 dt$$

$$= -\int_{-\infty}^{a} e^{2t}(\Omega - \Omega(a))^2 dt = -\int |u|^2 dt,$$

and a similar conclusion follows. This proves Proposition 4.2. QED.

## § 5. Perturbation theory. Part I

Theorem 5.1. Consider the eigenvalue problem

$$-\frac{d^2}{dt^2}\psi + m^2\psi + (\Omega(t) - \mu)^{-1}A(t)\psi = 0, \ t \in \mathbb{R}, \ Im \, \mu > 0, \ \psi \in L^2(\mathbb{R}), \psi \neq 0. \tag{5.1}$$

Let $\psi_0$ be an eigenfunction in (4.1) normalised so that $\|\psi_0\|_{L^2(\mathbb{R})} = 1$,

$$-\frac{d^2}{dt^2}\psi_0 + m_0{}^2\psi_0 + (\Omega(t) - \Omega(a))^{-1}A(t)\psi_0 = 0, \ t \in \mathbb{R}, \ \psi_0 \in L^2(\mathbb{R}), \ m_0 > 1. \tag{5.2}$$

Then, there exists an $\varepsilon > 0$ such that for any $h \in (0, \varepsilon)$ there is a solution $\psi = \psi_{m,\mu}$ of (5.1) with $m = m_0 - h$, such that $\mu = \mu(h) \to \Omega(a)$ as $h \to 0 +$. This solution $\psi_{m,\mu} \to \psi_0$ strongly in $L^2(\mathbb{R})$ as $h \to 0 +$.

Let $\psi_0$ be an eigenfunction in (4.1) normalized so that $\|\psi_0\|_{L^2(\mathbb{R})} = 1$,

$$-\frac{d^2}{dt^2}\psi_0 + m_0{}^2\psi_0 + (\Omega(t) - \Omega(b))^{-1}A(t)\psi_0 = 0, \ t \in \mathbb{R}, \ \psi_0 \in L^2(\mathbb{R}), \ m_0 > 1. \tag{5.3}$$

Then, likewise, there exists an $\varepsilon > 0$ such that for any $h \in (0, \varepsilon)$ there is a solution $\psi = \psi_{m,\mu}$ with $m = m_0 + h$, such that $\mu = \mu(h) \to \Omega(b)$ as $h \to 0 +$. This solution $\psi_{m,\mu} \to \psi_0$ strongly in $L^2(\mathbb{R})$ as $h \to 0 +$.

We give the following equivalent formulation of (5.1). Consider the eigenfunction

$$-\frac{d^2}{dt^2}\psi_0 + m_0{}^2\psi_0 + (\Omega(t) - \mu_0)^{-1}A(t)\psi_0 = 0, \ t \in \mathbb{R}, \ \psi_0 \in L^2(\mathbb{R}), \|\psi_0\|_{L^2(\mathbb{R})} = 1, \tag{5.4}$$

where $\mu_0 = \Omega(a)$ or $\mu_0 = \Omega(b)$. Notice $A(a) = A(b) = 0$. We use the following scheme of the perturbation theory.

Let $|h| < \varepsilon$,

$$\mu \in \mathcal{B} = \{z \in \mathbb{C} | \gamma < arg(z - \mu_0) < \pi - \gamma \}, \tag{5.5}$$

an angle, sufficiently close to $\pi$, a small constant $\gamma > 0$ will be specified later.

We replace the problem (5.1) by the following (see [FR])

$$-\frac{d^2}{dt^2}\psi + m^2\psi + (\Omega(t) - \mu)^{-1}A(t)\psi + cP_0\psi = c\psi_0, \ (\psi, \psi_0) = 1, \ \psi \in L^2(\mathbb{R}), \ Im\,\mu > 0. \tag{5.6}$$

Here $c > 0$ is a fixed real number, $P_0 : L^2(\mathbb{R}) \to \mathbb{R}\psi_0$

is a projector

$$P_0 y = (y, \psi_0)\psi_0, \quad \forall y \in L^2(\mathbb{R}). \tag{5.7}$$

It is clear that (5.6) is equivalent to (5.1). From (5.2), (5.6), (5.7)

$$-\frac{d^2}{dt^2}(\psi - \psi_0) + m_0^2(\psi - \psi_0) + cP_0(\psi - \psi_0) + (\Omega(t) - \mu_0)^{-1}A(t)(\psi - \psi_0)$$

$$+((\Omega(t) - \mu)^{-1} - (\Omega(t) - \mu_0)^{-1})A(t)\psi + (m^2 - m_0^2)\psi = 0, \quad (\psi, \psi_0) = 1. \tag{5.8}$$

For the proof of Theorem 5.1 that will be presented in § 7 we will need a few technical lemmas.

Lemma 5.1. For $\mu \in \mathcal{B}$, the operator

$$y \to ((\Omega(t) - \mu)^{-1} - (\Omega(t) - \mu_0)^{-1})A(t)y$$

is uniformly bounded in $\mathcal{L}(L^2(\mathbb{R}))$. As, $\mu \to \mu_0$, this operator converges to 0 strongly in $\mathcal{L}(L^2(\mathbb{R}))$.

Proof. We have,

$$|((\Omega(t) - \mu)^{-1} - (\Omega(t) - \mu_0)^{-1})A(t)|$$

$$= |\mu - \mu_0| \, |\Omega(t) - \mu|^{-1} \, |\Omega(t) - \mu_0|^{-1}|A(t)|$$

$$\leq C Im\, \mu \, |\Omega(t) - \mu|^{-1} \, |\Omega(t) - \mu_0|^{-1}|A(t)|$$

$$\leq C,$$

Since $|\mu - \mu_0| \leq C\, Im\, \mu$, $\forall \mu \in \mathcal{B}$. Let $d = a$ or $d = b$ so that $\Omega(d) = \mu_0$.

For any $\delta > 0$ and $|t - d| > \delta$

$$|\mu - \mu_0| \, |\Omega(t) - \mu|^{-1} \, |\Omega(t) - \mu_0|^{-1}|A(t)| \leq C\delta^{-2}(\sin\gamma)^{-2}\, |A(t)|\, Im\, \mu\, .$$

Lemma 5.1 follows from this inequality by Lebesgue dominated convergence theorem. QED.

Let

$$L_{m_0} = -\frac{d^2}{dt^2} + m_0^2\, ,$$

$$K_{m_0} = L_{m_0}^{-1}$$

in $L^2(\mathbb{R})$. The operator $K_{m_0} \in \mathcal{L}(L^2(\mathbb{R}))$ is an integral operator with kernel $K_{m_0}(\xi, \eta)$ (see (1.35)).

The problem (5.6) is equivalent to

$$\psi + \{K_{m_0} \circ ((\Omega - \mu_0)^{-1}A + cP_0) + K_{m_0} \circ ((\Omega - \mu)^{-1} - (\Omega - \mu_0)^{-1})A$$

$$+(m^2 - m_0^2)K_{m_0}\}\psi = cK_{m_0}\psi_0, \quad (\psi, \psi_0) = 1, \quad \psi \in L^2(\mathbb{R}), \quad Im\,\mu > 0. \tag{5.9}$$

The operator

$$id + K_{m_0} \circ ((\Omega - \mu_0)^{-1}A + cP_0) \tag{5.10}$$

is bounded and invertible in $\mathcal{L}(L^2(\mathbb{R}))$. It is also bounded and invertible in any Holder space $C^\sigma(\mathbb{R})$, $\sigma \in (0,1)$. We treat the third and the fourth term in the left side of (5.9) as perturbation.

Lemma 5.2. In $\mathcal{L}(L^2(\mathbb{R}))$

$$\| K_{m_0} \circ ((\Omega - \mu)^{-1} - (\Omega - \mu_0)^{-1})A \|_{\mathcal{L}(L^2(\mathbb{R}))} \to 0$$

as $\mu \in \mathfrak{B}$, $\mu \to \mu_0$.

Proof. Let $\chi \in C_0^\infty(\mathbb{R})$, $\chi|_{[d-1,d+1]} = 1$. We have,

$$K_{m_0} \circ ((\Omega - \mu)^{-1} - (\Omega - \mu_0)^{-1})A$$

$$= K_{m_0} \circ \chi((\Omega - \mu)^{-1} - (\Omega - \mu_0)^{-1})A + K_{m_0} \circ (1 - \chi)((\Omega - \mu)^{-1} - (\Omega - \mu_0)^{-1})A.$$

The second term is bounded in $\mathcal{L}(L^2(\mathbb{R}))$ and analytic in $\mu$ in some disc centered at $\mu_0$. As for the first term, $K_{m_0} \circ \chi \in \mathfrak{S}_\infty(L^2(\mathbb{R})) \subset \mathcal{L}(L^2(\mathbb{R}))$. Therefore, given any $\varepsilon > 0$, there exists an approximation by an operator of a finite rank, so that

$$\| K_{m_0} \circ \chi - \sum_{j=1}^{N}(\cdot, g_j)f_j \|_{\mathcal{L}(L^2(\mathbb{R}))} < \varepsilon.$$

Here $g_j, f_j \in L^2(\mathbb{R}), j = 1,2,3, \dots, N$.

Using Lemma 5.1 (applied to a complex conjugate function) with a finite number of vectors $g_j, j = 1, \dots, N$, we get

$$\| K_{m_0} \circ ((\Omega - \mu)^{-1} - (\Omega - \mu_0)^{-1})A \|_{\mathcal{L}(L^2(\mathbb{R}))} < 2\varepsilon$$

for $\mu \in \mathfrak{B}, |\mu - \mu_0| < \delta$, for a sufficiently small $\delta > 0$. This proves Lemma 5.2. QED.

Lemma 5.3. The operator $B_\mu$ acting in $C^\sigma(\mathbb{R})$, $\sigma \in (0,1)$

$$B_\mu : y \to (K_{m_0} \circ (\Omega - \mu)^{-1})y, \qquad \forall y \in C^\sigma(\mathbb{R})$$

is uniformly bounded in $\mathcal{L}(C^\sigma(\mathbb{R}))$ for $\mu \in \mathfrak{B}$. It converges strongly as $\mu \to \mu_0$, $\mu \in \mathfrak{B}$ to $B_{\mu_0+i0} \in \mathcal{L}(C^\sigma(\mathbb{R}))$, where

$$(B_{\mu_0+i0}\, y)(t) = \int K_{m_0}(t,\xi)(\Omega(\xi) - \mu_0 - i0)^{-1} y(\xi) d\xi \,.$$

Proof. Let $y \in C^\sigma(\mathbb{R})$. Choosing $\chi \in C_0^\infty(\mathbb{R})$ as in the proof of Lemma 5.2, we get

$$B_\mu y = K_{m_0} \circ \chi(\Omega - \mu)^{-1} y + K_{m_0} \circ (1-\chi)(\Omega - \mu)^{-1} y. \tag{5.11}$$

The second term in the right side of (5.11) is analytic with respect to $\mu \in \mathfrak{B} \cup B_r(\mu_0)$ with $r > 0$ sufficiently small. We turn to the first term in the right side of (5.11). We have (see (1.35))

$$[(K_{m_0} \circ \chi(\Omega - \mu)^{-1})y](t) = (2m_0)^{-1} \int e^{-m_0|t-\xi|} \chi(\xi)\, (\Omega(\xi) - \mu)^{-1} y(\xi) d\xi$$

$$= y(d)(2m_0)^{-1} \int e^{-m_0|t-\xi|} \chi(\xi)\, (\Omega(\xi) - \mu)^{-1} d\xi$$

$$+ (2m_0)^{-1} \int e^{-m_0|t-\xi|} \chi(\xi)\, (\Omega(\xi) - \mu)^{-1} (y(\xi) - y(d)) d\xi$$

$$= y(d)(2m_0)^{-1} I_1(t) + I_2(t) \,. \tag{5.12}$$

To estimate $I_2(t)$ we use the inequality

$$|\chi(\xi)||(\Omega(\xi) - \mu_0)|^{-1}|y(\xi) - y(d)| \le C\|y\|_{C^\sigma(\mathbb{R})} |\xi - d|^{\sigma - 1}, \qquad \forall \xi \in \mathrm{supp}\,\chi \,.$$

On the other hand, for $\mu \in \mathfrak{B}$, $|\Omega(\xi) - \mu| \ge \sin\gamma\, |\Omega(\xi) - \mu_0|$.

The function $t \to e^{-m_0|t|}$ belongs to $C^\sigma(\mathbb{R})$, $\forall \sigma \in (0,1)$. By Minkowski integral inequality,

$$\|I_2\|_{C^\sigma(\mathbb{R})} \le C\|y\|_{C^\sigma(\mathbb{R})}\,, \quad \forall \sigma \in (0,1)\,,$$

where the constant $C$ is uniform with respect to $\mu \in \mathfrak{B}$. We next estimate the integral $I_1$ in $C^\sigma(\mathbb{R})$.

$$I_1(t) = e^{-m_0|t-d|} \int \chi(\xi)\, (\Omega(\xi) - \mu)^{-1} d\xi$$

$$+ \int (e^{-m_0|t-\xi|} - e^{-m_0|t-d|})\, \chi(\xi)\, (\Omega(\xi) - \mu)^{-1} d\xi$$

$$= I_3(t) + I_4(t) \,. \tag{5.13}$$

We have

$$\|I_3\|_{C^\sigma(\mathbb{R})} \le C\,, \forall \sigma \in (0,1)$$

uniformly with respect to $\mu \in \mathfrak{B}$. As for the integral $I_4$,

$$I_4(t) = \int \int_d^\xi (m_0\, \mathrm{sign}\,(t-\rho)\, e^{-m_0|t-\rho|})\, \chi(\xi)\, (\Omega(\xi) - \mu)^{-1} d\rho d\xi \,.$$

Therefore,

$$\|I_4\|_{L^\infty(\mathbb{R})} \leq C \int |\chi(\xi)| \, |\xi - d| |\Omega(\xi) - \mu|^{-1} d\xi \leq C, \qquad (5.14)$$

uniformly with respect to $\mu \in \mathfrak{B}$. We estimate $|I_4(t_1) - I_4(t_2)|$ for $|t_1 - t_2| \leq 1$ by combining the following two estimates. First,

$$|e^{-m_0|t_1-\xi|} - e^{-m_0|t_1-d|} - e^{-m_0|t_2-\xi|} + e^{-m_0|t_2-d|}| \leq C|t_1 - t_2|.$$

Second,

$$|e^{-m_0|t_1-\xi|} - e^{-m_0|t_1-d|} - e^{-m_0|t_2-\xi|} + e^{-m_0|t_2-d|}| \leq C|\xi - d|.$$

Therefore, for any $\sigma \in (0,1)$

$$|e^{-m_0|t_1-\xi|} - e^{-m_0|t_1-d|} - e^{-m_0|t_2-\xi|} + e^{-m_0|t_2-d|}| \leq C|t_1 - t_2|^\sigma |\xi - d|^{1-\sigma}.$$

Hence,

$$|I_4(t_1) - I_4(t_2)| \leq C \, |t_1 - t_2|^\sigma \int |\chi(\xi)| \, |\xi - d|^{1-\sigma} |\Omega(\xi) - \mu|^{-1} d\xi$$

$$\leq C \, |t_1 - t_2|^\sigma. \qquad (5.15)$$

From (5.14), (5.15)

$$\|I_4\|_{C^\sigma(\mathbb{R})} \leq C, \forall \sigma \in (0,1),$$

where the constant $C$ is uniform with respect to $\mu \in \mathfrak{B}$. This completes the proof of the first statement in Lemma 5.3.

To prove strong convergence $B_\mu \to B_{\mu_0+i0}$ as $\mu \in \mathfrak{B}$, $\mu \to \mu_0$ we must pass to the limit in the first term in (5.11). Indeed, convergence (even in the uniform norm in $\mathcal{L}(C^\sigma(\mathbb{R}))$) is obvious in the second term because of analyticity in $\mu$ at $\mu_0$. We have

$$[(K_{m_0} \circ \chi(\Omega - \mu)^{-1})y](t) = y(d)(2m_0)^{-1}(I_3(t) + I_4(t)) + I_2(t).$$

The strong convergence of $I_2$ to the limit

$$(2m_0)^{-1} \int e^{-m_0|t-\xi|} \chi(\xi) \, (\Omega(\xi) - \mu_0)^{-1} (y(\xi) - y(d)) d\xi$$

$$= (2m_0)^{-1} \int e^{-m_0|t-\xi|} \chi(\xi) \, (\Omega(\xi) - \mu_0 - i0)^{-1} (y(\xi) - y(d)) d\xi$$

follows from the Lebesgue dominated convergence theorem. As for $I_3$ the statement holds since $\chi \in C_0^\infty(\mathbb{R})$ and $(\Omega - \mu)^{-1} \to (\Omega - \mu_0 - i0)^{-1}$ in $S'(\mathbb{R})$ as $\mu \in \mathfrak{B}$, $\mu \to \mu_0$. For $I_4$ the strong convergence in $L^\infty(\mathbb{R})$ and in the (semi-) norm involving

$$\sup_{0<|t_1-t_2|\leq 1} |t_1 - t_2|^{-\sigma}|I_4(t_1) - I_4(t_2)|$$

follows from Lebesgue theorem. This completes the proof of Lemma 5.3. QED.

# § 6. Perturbation theory. Part II

Lemma 6.1. Let $C^\sigma(\mathbb{R})$, $\sigma \in (0,1)$ be the Holder space with exponent $\sigma$. Let $m_0 \in (0,\infty)$. Define for $\mu \in H = \{z \in \mathbb{C} | \, Im \, z > 0\}$ the operator

$$B_\mu : y \to (K_{m_0} \circ (\Omega - \mu)^{-1})y, \qquad \forall y \in C^\sigma(\mathbb{R}).$$

As $\mu \to \mu_0 = \Omega(\nu_0)$ in $H$ for some fixed $\nu_0 \in \mathbb{R}$, the operator-valued function $B_\mu$ converges in the uniform sense in $\mathcal{L}(C^\sigma(\mathbb{R}))$ to $B_{\mu_0+i0} \in \mathcal{L}(C^\sigma(\mathbb{R}))$, where

$$(B_{\mu_0+i0} \, y)(t) = \int K_{m_0}(t,\xi)(\Omega(\xi) - \mu_0 - i0)^{-1}y(\xi)d\xi, \quad \forall y \in C^\sigma(\mathbb{R}).$$

Proof. We choose $\chi \in C_0^\infty(\mathbb{R})$, $\chi|_{[\nu_0-1,\nu_0+1]} \equiv 1$, and break $B_\mu$ as follows

$$B_\mu y = K_{m_0} \circ \chi(\Omega - \mu)^{-1}y + K_{m_0} \circ (1-\chi)(\Omega - \mu)^{-1}y, \qquad \forall y \in C^\sigma(\mathbb{R}). \tag{6.1}$$

The second term in the right side of (6.1) is analytic in $\mu \in H \cup B_r(\mu_0)$ for $r > 0$ sufficiently small. We turn to the first term

$$\left(K_{m_0} \circ \chi(\Omega - \mu)^{-1}y\right)(t) =$$

$$= (2m_0)^{-1} \int e^{-m_0|t-\xi|}\chi(\xi)(\Omega(\xi) - \mu)^{-1}y(\xi)d\xi$$

$$= y(\nu)(2m_0)^{-1} \int e^{-m_0|t-\xi|}\chi(\xi)(\Omega(\xi) - \mu)^{-1}d\xi$$

$$+ (2m_0)^{-1} \int e^{-m_0|t-\xi|}\chi(\xi)(\Omega(\xi) - \mu)^{-1}(y(\xi) - y(\nu))d\xi$$

$$= y(\nu)(2m_0)^{-1} I_1(t) + I_2(t), \tag{6.2}$$

where $\mu = \Omega(\nu) + i\kappa$, $\kappa > 0$. For the integral $I_2$ we use the inequality

$$|\chi(\xi)||(\Omega(\xi) - \mu)|^{-1}|y(\xi) - y(\nu)| \leq C\|y\|_{C^\sigma(\mathbb{R})} |\xi - \nu|^{\sigma-1}, \qquad \forall \xi \in supp\, \chi.$$

The function $t \to e^{-m_0|t|}$ belongs to $C^\sigma(\mathbb{R})$, $\forall \sigma \in (0,1)$. By Minkowski integral inequality,

$$\|I_2\|_{C^\sigma(\mathbb{R})} \leq C\|y\|_{C^\sigma(\mathbb{R})}, \quad \forall \sigma \in (0,1) \tag{6.3}$$

where the constant $C$ is uniform with respect to $\mu \in B_r(\mu_0) \cap H$.

We next estimate the integral $I_1$ in $C^\sigma(\mathbb{R})$.

$$I_1(t) = e^{-m_0|t-\nu|} \int \chi(\xi)\, (\Omega(\xi) - \mu)^{-1} d\xi$$

$$+ \int (e^{-m_0|t-\xi|} - e^{-m_0|t-\nu|})\, \chi(\xi)\, (\Omega(\xi) - \mu)^{-1} d\xi$$

$$= I_3(t) + I_4(t). \tag{6.4}$$

We have $\|I_3\|_{C^\sigma(\mathbb{R})} \leq C$ uniformly in $\mu \in B_r(\mu_0) \cap H$.

To estimate $\|I_4\|_{L^\infty(\mathbb{R})}$ we use the identity

$$I_4(t) = \int \int_\nu^\xi (m_0\, sign\, (t - \rho)\, e^{-m_0|t-\rho|})\, \chi(\xi)\, (\Omega(\xi) - \mu)^{-1} d\rho d\xi.$$

Therefore,

$$\|I_4\|_{L^\infty(\mathbb{R})} \leq C \int |\chi(\xi)||(\Omega(\xi) - \mu)|^{-1}|\xi - \nu|\, d\xi$$

$$\leq C \int |\chi(\xi)||(\Omega(\xi) - \mu)|^{-1}|\xi - \nu|\, d\xi \leq C, \tag{6.5}$$

where the constant $C$ is uniform with respect to $\mu \in B_r(\mu_0) \cap H$.

We estimate $|I_4(t_1) - I_4(t_2)|$ for $|t_1 - t_2| \leq 1$ by combining the following two estimates. First,

$$|e^{-m_0|t_1-\xi|} - e^{-m_0|t_1-\nu|} - e^{-m_0|t_2-\xi|} + e^{-m_0|t_2-\nu|}| \leq C|t_1 - t_2|.$$

Second,

$$|e^{-m_0|t_1-\xi|} - e^{-m_0|t_1-\nu|} - e^{-m_0|t_2-\xi|} + e^{-m_0|t_2-\nu|}| \leq C|\xi - \nu|.$$

Therefore, for any $\sigma \in (0,1)$

$$|e^{-m_0|t_1-\xi|} - e^{-m_0|t_1-\nu|} - e^{-m_0|t_2-\xi|} + e^{-m_0|t_2-\nu|}| \leq C|t_1 - t_2|^\sigma |\xi - \nu|^{1-\sigma}.$$

Thus, for $|t_1 - t_2| \leq 1$

$$|I_4(t_1) - I_4(t_2)| \leq C \ |t_1 - t_2|^\sigma \int |\xi - v|^{1-\sigma} |\chi(\xi)| |(\Omega(\xi) - \mu)|^{-1} d\xi \leq C \ |t_1 - t_2|^\sigma \tag{6.6}$$

uniformly in $\mu \in B_r(\mu_0) \cap H$ with $r > 0$ small enough.

From (6.5), (6.6) we have $\|I_4\|_{C^\sigma(\mathbb{R})} \leq C$ uniformly in $\mu \in B_r(\mu_0) \cap H$.

Therefore, the family of operators $B_\mu \in \mathcal{L}(C^\sigma(\mathbb{R}))$ is uniformly bounded for $\mu \in B_r(\mu_0) \cap H$ with $r > 0$ small enough. It remains to show convergence $B_\mu \to B_{\mu_0 + i0}$ in the uniform sense in $\mathcal{L}(C^\sigma(\mathbb{R}))$.

Since passing to the limit in the second term in the right side of (6.1) is obvious, we concentrate on the first term $y(v)(2m_0)^{-1} I_1(t) + I_2(t)$. To pass to the limit in $I_2$ we notice (making the dependence on $\mu$ explicit)

$$I_2(t, \mu) - \lim_{\mu \to \mu_0 + i0} I_2(t, \mu)$$

$$= (2m_0)^{-1} \int e^{-m_0 |t - \xi|} \chi(\xi) \left[ (\Omega(\xi) - \mu)^{-1} (y(\xi) - y(v)) - (\Omega(\xi) - \mu_0)^{-1} (y(\xi) - y(v_0)) \right] d\xi.$$

Let $\|y\|_{C^\sigma(\mathbb{R})} = 1$. We claim

$$\|\chi(\xi)[(\Omega(\xi) - \mu)^{-1}(y(\xi) - y(v)) - (\Omega(\xi) - \mu_0)^{-1}(y(\xi) - y(v_0))]\|_{L^1(\mathbb{R})} \to 0$$

as $\mu \to \mu_0$ in $B_r(\mu_0) \cap H$ uniformly with respect to $y$. Here $r > 0$ is small enough. Indeed,

$$g(\xi) \equiv \chi(\xi)[(\Omega(\xi) - \mu)^{-1}(y(\xi) - y(v)) - (\Omega(\xi) - \mu_0)^{-1}(y(\xi) - y(v_0))]$$

$$= \chi(\xi) \left( [(\Omega(\xi) - \mu)^{-1} - (\Omega(\xi) - \mu_0)^{-1}](y(\xi) - y(v)) - (\Omega(\xi) - \mu_0)^{-1}(y(v) - y(v_0)) \right).$$

Choosing a small $\delta > 0$ we have

$$\int_{|\xi - v_0| > \delta} |g(\xi)| d\xi \leq C(\delta)(|\mu - \mu_0| + |\mu - \mu_0|^\sigma). \tag{6.6}$$

Also,

$$\int_{|\xi - v_0| \leq \delta} |g(\xi)| d\xi \leq C \delta^\sigma$$

as $|v - v_0| < \frac{1}{2} \delta$. The right side of (6.6) can be made less than $\delta^\sigma$ as $\mu \to \mu_0$ in $B_r(\mu_0) \cap H$, uniformly in $y$.

By Minkowski integral inequality

$$\|I_2(t,\mu) - \lim_{\mu \to \mu_0 + i0} I_2(t,\mu)\|_{\mathcal{L}(C^\sigma(\mathbb{R}))} \to 0$$

as $\mu \to \mu_0$ in $B_r(\mu_0) \cap H$.

To prove convergence in the first term in the right side of (6.2) it is sufficient to prove strong convergence of $I_1(t,\mu)$ in $C^\sigma(\mathbb{R})$ as $\mu \to \mu_0$ in $B_r(\mu_0) \cap H$. Indeed,

$$|y(\nu) - y(\nu_0)| \leq C\|y\|_{C^\sigma(\mathbb{R})}|\mu - \mu_0|^\sigma.$$

This reduces the statement to strong convergence of $I_3(t,\mu)$ and $I_4(t,\mu)$. As for $I_3(t,\mu)$ we have

$$\int \chi(\xi)(\Omega(\xi) - \mu)^{-1}d\xi \to \int \chi(\xi)(\Omega(\xi) - \mu_0 - i0)^{-1}d\xi,$$

since $\chi \in C_0^\infty(\mathbb{R})$ and $(\Omega - \mu)^{-1} \to (\Omega - \mu_0 - i0)^{-1}$ in $\mathcal{D}'(\mathbb{R})$. Also it is easy to see that

$$e^{-m_0|.-\nu|} \to e^{-m_0|.-\nu_0|}$$

strongly in $C^\sigma(\mathbb{R})$ as $\nu \to \nu_0$.

As for the term $I_4$,

$$\int (e^{-m_0|t-\xi|} - e^{-m_0|t-\nu|})\chi(\xi)(\Omega(\xi) - \mu)^{-1} - (e^{-m_0|t-\xi|} - e^{-m_0|t-\nu_0|})\chi(\xi)(\Omega(\xi) - \mu_0)^{-1} d\xi$$

$$= \int |\xi - \nu_0|^{\sigma-1}\left(e^{-m_0|t-\xi|} - e^{-m_0|t-\nu_0|}\right)\chi(\xi)$$

$$(|\xi - \nu|^{1-\sigma}(\Omega(\xi) - \mu)^{-1} - |\xi - \nu_0|^{1-\sigma}(\Omega(\xi) - \mu_0)^{-1})d\xi$$

$$+ \int \{|\xi - \nu|^{\sigma-1}\left(e^{-m_0|t-\xi|} - e^{-m_0|t-\nu|}\right) - |\xi - \nu_0|^{\sigma-1}\left(e^{-m_0|t-\xi|} - e^{-m_0|t-\nu_0|}\right)\}$$

$$\chi(\xi)|\xi - \nu|^{1-\sigma}(\Omega(\xi) - \mu)^{-1}d\xi. \tag{6.7}$$

In the first integral in the right side of (6.7)

$$\|\{t \to |\xi - \nu_0|^{\sigma-1}\left(e^{-m_0|t-\xi|} - e^{-m_0|t-\nu_0|}\right)\}\|_{C^\sigma(\mathbb{R})} \leq C,$$

where the constant $C$ is independent of $\xi$ and $\nu_0$. Also, if

$$g(\xi) = \chi(\xi)(|\xi - \nu|^{1-\sigma}(\Omega(\xi) - \mu)^{-1} - |\xi - \nu_0|^{1-\sigma}(\Omega(\xi) - \mu_0)^{-1}),$$

than $\|g\|_{L^1(\mathbb{R})} \to 0$ as $\mu \to \mu_0$ in $B_r(\mu_0) \cap H$. Indeed, this follows from the Lebesgue dominated convergence theorem. Therefore, by Minkowski integral inequality the first integral in the right side of (6.7) as a function of $t$ converges to zero strongly in $C^\sigma(\mathbb{R})$ as $\mu \to \mu_0$ in $B_r(\mu_0) \cap H$.

In the second integral in the right side of (6.7) by the same reasoning we can replace
$|\xi - v|^{1-\sigma}(\Omega(\xi) - \mu)^{-1}$ by $|\xi - v_0|^{1-\sigma}(\Omega(\xi) - \mu_0)^{-1}$. We split the resulting integral as follows:

$$\int \{|\xi - v|^{\sigma-1}(e^{-m_0|t-\xi|} - e^{-m_0|t-v|}) - |\xi - v_0|^{\sigma-1}(e^{-m_0|t-\xi|} - e^{-m_0|t-v_0|})\}$$

$$\chi(\xi)|\xi - v_0|^{1-\sigma}(\Omega(\xi) - \mu_0)^{-1}d\xi$$

$$= \int_{|\xi-v_0|>\delta} \{|\xi - v|^{\sigma-1}(e^{-m_0|t-\xi|} - e^{-m_0|t-v|}) - |\xi - v_0|^{\sigma-1}(e^{-m_0|t-\xi|} - e^{-m_0|t-v_0|})\}\chi(\xi)$$

$$|\xi - v_0|^{1-\sigma}(\Omega(\xi) - \mu_0)^{-1}d\xi$$

$$+ \int_{|\xi-v_0|\leq\delta} \{|\xi - v|^{\sigma-1}(e^{-m_0|t-\xi|} - e^{-m_0|t-v|}) - |\xi - v_0|^{\sigma-1}(e^{-m_0|t-\xi|} - e^{-m_0|t-v_0|})\}\chi(\xi)$$

$$|\xi - v_0|^{1-\sigma}(\Omega(\xi) - \mu_0)^{-1}d\xi$$

$$= I_5 + I_6 .$$

For $I_6$ we have $\|I_6\|_{C^\sigma(\mathbb{R})} \leq C\delta^{1-\sigma}$. In the integral $I_5$

$$|\xi - v|^{\sigma-1}(e^{-m_0|t-\xi|} - e^{-m_0|t-v|}) - |\xi - v_0|^{\sigma-1}(e^{-m_0|t-\xi|} - e^{-m_0|t-v_0|})$$

$$= |\xi - v_0|^{\sigma-1}(e^{-m_0|t-v_0|} - e^{-m_0|t-v|})$$

$$+(|\xi - v|^{\sigma-1} - |\xi - v_0|^{\sigma-1})(e^{-m_0|t-\xi|} - e^{-m_0|t-v|}) .$$

Therefore,

$$\|I_5\|_{C^\sigma(\mathbb{R})} \leq C\delta^{-(1-\sigma)}|v - v_0|^{1-\sigma} + C\delta^{-(2-\sigma)}|v - v_0|,$$

provided $|v - v_0| < \frac{\delta}{2}$. This does to 0 as $v \to v_0$. Since $\delta > 0$ is arbitrary small, we have

$$\|I_5\|_{C^\sigma(\mathbb{R})} \to 0$$

as $\mu \to \mu_0$ in $B_r(\mu_0) \cap H$. This completes the proof of Lemma 6.1, QED.

## §7. Proof of Theorem 5.1

By Lemma 5.2 we can now define for $\mu \in \mathfrak{B}$; $|m - m_0|, |\mu - \mu_0|$ small enough,

$\psi = \psi_{m,\mu} \in L^2(\mathbb{R})$ as a unique solution of the equation (5.9) disregarding the normalization condition $(\psi_{m,\mu}, \psi_0) = 1$. This condition becomes the characteristic equation. Evidently, $\psi = \psi_{m,\mu}$ is analytic in $\mu \in \mathfrak{B} \cap B_r(\mu_0)$ for a sufficiently small $r > 0$. By this construction

$$\psi_{m,\mu} \to \psi_0$$

strongly in $L^2(\mathbb{R})$ as $m \to m_0$, $\mu \to \mu_0$; $\mu \in \mathfrak{B}$. Notice $\psi = \psi_{m,\mu}$ satisfies (5.8) [without normalization condition $(\psi_{m,\mu}, \psi_0) = 1$.]

A similar statement holds in $C^\sigma(\mathbb{R})$ for any $\sigma \in (0,1)$: $\psi_{m,\mu} \to \psi_0$ strongly in $C^\sigma(\mathbb{R})$ as

$$m \to m_0, \qquad \mu \to \mu_0; \quad \mu \in H.$$

Indeed,

$$K_{m_0} \circ ((\Omega - \mu)^{-1} - (\Omega - \mu_0)^{-1})A$$

$$= (\mu - \mu_0) K_{m_0} \circ (\Omega - \mu)^{-1} \circ (\Omega - \mu_0)^{-1}A \to 0$$

in the uniform sense in $\mathcal{L}(C^\sigma(\mathbb{R}))$ as $\mu \to \mu_0, \mu \in H$, by Lemma 6.1. We note that $(\Omega - \mu_0)^{-1}A \in \mathcal{L}(C^\sigma(\mathbb{R}))$ since $\mu_0 = \Omega(a)$ or $\mu_0 = \Omega(b)$ in the statement of Theorem 5.1.

Indeed, from Lemma 6.1

$$\|K_{m_0} \circ (\Omega - \mu)^{-1} \circ (\Omega - \mu_0)^{-1}A\|_{\mathcal{L}(C^\sigma(\mathbb{R}))} = O(1)$$

as $\mu \to \mu_0, \mu \in H$.

Now from the Neumann series expansion in (5.9)

$$\psi_{m,\mu} = \sum_{n=0}^{\infty}(-1)^n \{(id + K_{m_0}((\Omega - \mu_0)^{-1}A + cP_0))^{-1}$$

$$[K_{m_0}((\Omega - \mu)^{-1} - (\Omega - \mu_0)^{-1})A + (m^2 - m_0^2)K_{m_0}]\}^n \psi_0.$$

Indeed, by construction,

$$(id + K_{m_0}((\Omega - \mu_0)^{-1}A + cP_0))^{-1} cK_{m_0}\psi_0 = \psi_0. \tag{7.1}$$

Therefore, keeping two terms in the Neumann series, we get the following inequality:

$\|\psi_{m,\mu} - \psi_0 - (id + K_{m_0}((\Omega - \mu_0)^{-1}A + cP_0))^{-1}$

$\qquad [K_{m_0}((\Omega - \mu)^{-1} - (\Omega - \mu_0)^{-1})A + (m^2 - m_0^2)K_{m_0}]\psi_0\|_{C^\sigma(\mathbb{R})} \leq C(|\mu - \mu_0| + |m - m_0|)^2.$

Using Lemma 6.1 again, we arrive at

$$\|\psi_{m,\mu} - \psi_0 - (id + K_{m_0}((\Omega - \mu_0)^{-1}A + cP_0))^{-1}$$

$$\left[(\mu - \mu_0)K_{m_0}(\Omega - \mu_0 - i0)^{-1}(\Omega - \mu_0)^{-1}A + 2m_0(m - m_0)K_{m_0}\right]\psi_0 \|_{C^\sigma(\mathbb{R})}$$

$$\leq F(|\mu - \mu_0| + |m - m_0|),$$

where $\mu \in H, F(s) = o(s)$ as $s \to 0+$.

Therefore, since $\|\psi_0\|_{L^2(\mathbb{R})} = 1$ (see (5.4)),

$$(\psi_{m,\mu}, \psi_0) - 1 - ((id + K_{m_0}((\Omega - \mu_0)^{-1}A + cP_0))^{-1}$$

$$\left[(\mu - \mu_0)K_{m_0}(\Omega - \mu_0 - i0)^{-1}(\Omega - \mu_0)^{-1}A + 2m_0(m - m_0)K_{m_0}\right]\psi_0, \psi_0) = o(|\mu - \mu_0| + |m - m_0|)$$

as $\mu \in H$.

From (5.4), (7.1) and using

$$(id + ((\Omega - \mu_0)^{-1}A + cP_0)K_{m_0})^{-1}\psi_0 = c^{-1}L_{m_0}\psi_0,$$

we get

$$(\psi_{m,\mu}, \psi_0) - 1 - 2m_0(m - m_0)c^{-1} - (\mu - \mu_0)c^{-1}(K_{m_0}(\Omega - \mu_0 - i0)^{-1}[(\Omega - \mu_0)^{-1}A]\psi_0, L_{m_0}\psi_0)$$

$$= o(|\mu - \mu_0| + |m - m_0|), \quad \mu \in H.$$

Thus,

$$c((\psi_{m,\mu}, \psi_0) - 1) - 2m_0(m - m_0) - (\mu - \mu_0)((\Omega - \mu_0 - i0)^{-1}[(\Omega - \mu_0)^{-1}A]\psi_0, \psi_0)$$

$$= o(|\mu - \mu_0| + |m - m_0|), \quad \mu \in H. \quad (7.2)$$

Let $m = m_0 + h$, $\mu = \mu_0 - (2m_0 + z)h \, ((\Omega - \mu_0 - i0)^{-1}[(\Omega - \mu_0)^{-1}A]\psi_0, \psi_0)^{-1}$,

where $|z| = r$ is small enough. The imaginary part of the denominator does not vanish.

Indeed, for $d = a$ or $d = b$

$$Im \, ((\Omega - \mu_0 - i0)^{-1}[(\Omega - \mu_0)^{-1}A]\psi_0, \psi_0) = \pi \frac{1}{|\Omega'(d)|} [(\Omega - \mu_0)^{-1}A](d)|\psi_0(d)|^2. \quad (7.3)$$

But $\psi_0(d) \neq 0$, otherwise $\psi_0 \equiv 0$ on $\mathbb{R}$. Indeed, see the proof of Proposition (4.1). Therefore,

$$Im \, ((\Omega - \Omega(b) - i0)^{-1}[(\Omega - \mu_0)^{-1}A](b)\psi_0, \psi_0) > 0,$$

$$Im \, ((\Omega - \Omega(a) - i0)^{-1}[(\Omega - \mu_0)^{-1}A](a)\psi_0, \psi_0) < 0.$$

To have $\mu \in H$ we choose

$$h < 0 \, for \, d = a;$$

$$h > 0 \ for \ d = b.$$

If we choose $r \in (0, \ 2m_0)$, the disc, centered at

$$\mu_0 - 2m_0 h \, ((\Omega - \mu_0 - i0)^{-1}, [(\Omega - \mu_0)^{-1} A] |\psi_0|^2)^{-1}$$

of radius

$$rh |((\Omega - \mu_0 - i0)^{-1}, [(\Omega - \mu_0)^{-1} A] |\psi_0|^2)|^{-1}$$

belongs entirely to the upper half-plane $H$. The function $z \to (\psi_{m,\mu}, \psi_0) - 1$ that is analytic on this disc for a fixed $h$, $|h| < \varepsilon$ equals $-c^{-1} zh + o(|h|)$ on the boundary of the disc. According to Rouchet's theorem there is exactly one zero of this function in the interior of this disc. This completes the proof of Theorem 5.1. QED.

## §8. Bottom of the spectrum

We recast the description of the class $\mathcal{C}$ in terms of the vorticity $G$:

$(i)$ $G \in C^\infty(\mathbb{R})$;

$(ii)$ $\int_{-\infty}^{t} e^{2\tau} G'(\tau) d\tau < 0, \forall t \in \mathbb{R}$;

$(iii)$ $G(t) = 2 \, \Omega(-\infty) - 4 c_0 e^{2t}$ for some $c_0 > 0$, $\Omega(-\infty) > 0$ and $t \leq \log M_1$;

$(iv)$ $G(t) = e^{-\alpha t}$ for $t \geq \log M > \log M_1$ ;

$(v)$ $A(t) = G'(t)$ has exactly 2 zeroes, $a$ and $b$, $\log M_1 < a < b < \log M$,

$\quad G''(a) > 0, \ G''(b) < 0$.

Remark 8.1. $\Omega(t) = \int_{-\infty}^{t} e^{-2(t-\tau)} G(\tau) d\tau$ is strictly decreasing on $\mathbb{R}$. Indeed,

$\Omega'(t) = \int_{-\infty}^{t} e^{-2(t-\tau)} G'(\tau) d\tau < 0, \forall t \in \mathbb{R}$ .

Remark 8.2. Let $a$ and $b$, $a < b$ be fixed, the set of functions $A(t) = G'(t)$ with zeroes $a$ and $b$ as above, corresponding to the class $\mathcal{C}$ is convex in $C^\infty(\mathbb{R})$.

Proposition 8.1. Let $N > 0$ be fixed. There exists a function $\Omega \in \mathcal{C}$, so that the infimum of the spectrum in $L^2(\mathbb{R})$ of the operator

$$L \equiv -\frac{d^2}{dt^2} + \left(\Omega(t) - \Omega(a)\right)^{-1} A(t)$$

is $\leq -N^2$.

Remark 8.3. The bottom of the spectrum of $L$ in $L^2(\mathbb{R})$ is a simple eigenvalue.

Proof of Proposition 8.1. We fix $a = 0$, $b = 1$. We construct $A(t)$ in several steps.

1. Let $B > 0$ be a large parameter. We define $A(t)$ for $t \in [-10^{-2} B^{-\frac{1}{2}}, 0]$ as follows:

$$A(t) = (4 + B)t + Bt^2, \ t \in [-10^{-2} B^{-\frac{1}{2}}, 0].$$

2. Extend $A(t)$ to $t$ $(-\infty, 0]$ so that $A(t) = -8c_0 e^{2t}$ for some $c_0 > 0$ and $t \in (-\infty, \log M_1]$,

$$A \in C^\infty((-\infty; 0); (-\infty, 0)); \ \int_{-\infty}^{0} e^{2\tau} A(\tau) d\tau = -1.$$

Notice that $\int_{-10^{-2} B^{-\frac{1}{2}}}^{0} A(\tau) d\tau = -\frac{1}{2}(4 + B)10^{-4} B^{-1} + \frac{1}{3} 10^{-6} B^{-\frac{1}{2}} \in (-10^{-4}, 0)$

for $B > 0$ sufficiently large.

3. Extend $A(t)$ to $\mathbb{R}$ so that $A \in C^\infty(\mathbb{R}; \mathbb{R})$, $A: (0,1) \to (0, e^{-2})$, $A: (1, \infty) \to (-\infty, 0)$;

$$A(t) = -\alpha e^{-\alpha t}, \forall t \in (\log M, \infty).$$

4. Define

$$G(t) = -\int_t^\infty A(\tau) d\tau, \ \forall t \in \mathbb{R}.$$

Then,

$$G(t) = e^{-\alpha t}, \quad \forall t \in (\log M, \infty),$$

$$G(t) = 2\Omega(-\infty) - 4c_0 e^{2t}, \quad \forall t \in (-\infty, \log M_1);$$

for some $\Omega(-\infty)$;

$$2\Omega(-\infty) = -\int_{-\infty}^{\infty} A(\tau) d\tau.$$

Since $\int_{-\infty}^{0} A(\tau)d\tau < -1$, $\int_{0}^{1} A(\tau)d\tau \in (0, e^{-2})$, $\int_{1}^{\infty} A(\tau)d\tau < 0$, it follows that $\Omega(-\infty) > 0$. We have

$$\Omega(-\infty) > 0.$$

Define now

$\Omega(t) = \int_{-\infty}^{t} e^{-2(t-\tau)} G(\tau) \, d\tau$, $\forall t \in \mathbb{R}$.

Then,

$$\Omega'(t) = \int_{-\infty}^{t} e^{-2(t-\tau)} G'(\tau) \, d\tau, \forall t \in \mathbb{R}.$$

Therefore, for $t \in (0,1]$

$$\Omega'(t) < e^{-2t}(-1 + \frac{1}{2}e^{2t-2}) < 0.$$

Thus $\Omega'(t) < 0$, $\forall t \in \mathbb{R}$. Also, $\Omega' + 2\Omega = G$, therefore,

$$A = \Omega'' + 2\Omega'.$$

On the interval $t \in [-10^{-2}B^{-\frac{1}{2}}, 0]$

$$\Omega(t) = \Omega(0) - t + t^2 + \frac{1}{6}Bt^3 + C(e^{-2t} - 1)$$

With some $C \in \mathbb{R}$. Indeed, for this $\Omega(t)$

$$\Omega'' + 2\Omega' = -2 + 2 + 4t + Bt + Bt^2 = (4 + B)t + Bt^2 = A(t).$$

Since $\Omega'(0) = -1$ by construction, we have $C = 0$. Thus, on the interval $\left[-10^{-2}B^{-\frac{1}{2}}, 0\right]$ we have

$$(\Omega(t) - \Omega(a))^{-1} A(t) = \left(-t + t^2 + \frac{1}{6}Bt^3\right)^{-1} ((4 + B)t + Bt^2)$$

$$= \left(-1 + t + \frac{1}{6}Bt^2\right)^{-1} (4 + B + Bt), \quad \forall t \in \left[-10^{-2}B^{-\frac{1}{2}}, 0\right].$$

We use the following test function in a standard calculus of variations argument:

$$\eta(t) = \begin{cases} \cos\left(\frac{\pi}{2}t\right), & t \in [-1,1] \\ 0, & t \in \mathbb{R}\setminus[-1,1] \end{cases}.$$

Obviously, $\eta \in W^{1,2}(\mathbb{R})$. We have

$$E(\eta) = \int_{\mathbb{R}} |\frac{d\eta}{dt}|^2 + (\Omega(t) - \Omega(a))^{-1} A(t)|\eta|^2 \, dt$$

$$= \frac{\pi^2}{4} + \int_{-1}^{1} (\Omega(t) - \Omega(a))^{-1} A(t)|\eta|^2 \, dt \,.$$

Notice that $(\Omega(t) - \Omega(a))^{-1} A(t) < 0$, $\forall t \in (-\infty, 1)$. Therefore,

$$\int_{-1}^{1} (\Omega(t) - \Omega(a))^{-1} A(t)|\eta|^2 \, dt < \int_{-10^{-2}B^{-\frac{1}{2}}}^{0} (\Omega(t) - \Omega(a))^{-1} A(t) \cos^2\left(\frac{\pi}{2} t\right) dt$$

$$< -10^{-2} B^{-\frac{1}{2}} \frac{1}{2}\left(1 + \frac{3}{2B}\right)^{-1} \left(4 + B - B^{\frac{1}{2}}\right) < -CB^{\frac{1}{2}}$$

for an appropriate $C > 0$. But $\|\eta\|_{L^2(\mathbb{R})} = 1$. Thus,

$$E(\eta) \leq -CB^{\frac{1}{2}} + \frac{\pi^2}{4} \leq -N^2 \, \|\eta\|_{L^2(\mathbb{R})}^2$$

for sufficiently large $B > 0$. This completes the proof of Proposition 8.1. QED.

## §9. $\mathcal{U}_m = \emptyset$ for large $m > 0$

Proposition 9.1. There exists $m_* > 0$ such that $\mathcal{U}_m = \emptyset$ for any $m \geq m_*$.

Proof. Assume the contrary.

The compliment $\mathbb{R} \setminus [0, \Omega(-\infty)]$ cannot contain a limit of a sequence $\mu_j \in \mathcal{U}_{m_j}$ as $m_j \to \infty$; $j = 1,2,3,\ldots$. Such a limit cannot be $\infty$ either. Indeed, if $\psi_j$

$$-\frac{d^2}{dt^2}\psi_j + m_j^2 \psi_j + (\Omega(t) - \mu_j)^{-1} A(t)\psi_j = 0, \ t \in \mathbb{R}, \ \text{Im}\, \mu_j > 0, \ \psi_j \in L^2(\mathbb{R}), \psi_j \neq 0, \qquad (9.1)$$

implies for a $\delta > 0$, $\delta \leq dist(\mu_j, [0, \Omega(-\infty)])$

$$\int |\frac{d}{dt}\psi_j|^2 + m_j^2 |\psi_j|^2 \, dt \leq \delta^{-1} C \int |\psi_j|^2 \, dt.$$

This leads to a contradiction for sufficiently large $j$, as $m_j \to \infty$.

Therefore, there exists a $\mu_0 \in [0, \Omega(-\infty)]$ such that

$$\mu_0 = \lim_{j \to \infty} \mu_j, \ \mu_j \in \mathcal{U}_{m_j}, \ \forall j = 1,2,3,\ldots; \ m_j \to \infty \text{ as } j \to \infty.$$

We first rule out the case $\mu_0 = 0$. We normalize $\psi_j \in L^2(\mathbb{R})$ so that

$$\psi_j(t) \sim e^{m_j t} \text{ as } t \to -\infty.$$

Selecting $any\ N \geq b$ and integrating on $(-\infty, \infty)$ we get

$$-\int_{-\infty}^{a} \left(|\Omega(t) - Re\mu_j|^2 + |Im\ \mu_j|^2\right)^{-1} A(t)|\psi_j(t)|^2\ dt$$

$$-\int_{b}^{N} \left(|\Omega(t) - Re\mu_j|^2 + |Im\ \mu_j|^2\right)^{-1} A(t)|\psi_j(t)|^2\ dt$$

$$\leq \int_{a}^{b} \left(|\Omega(t) - Re\mu_j|^2 + |Im\ \mu_j|^2\right)^{-1} A(t)|\psi_j(t)|^2\ dt. \tag{9.2}$$

Notice, the left side of this inequality is positive. For a fixed $N \geq b$ let $j$ be large enough, so that $|\mu_j| < \frac{1}{2}\Omega(N)$. From Lemma 2.1

$$\psi_j(t) = e^{m_j t}\left(1 + z_j(t)\right),$$

where

$$|z_j(t)| \leq \exp\left\{\frac{1}{2m_j}\int_{-\infty}^{t} |\Omega(\tau) - \mu_j|^{-1} |A(\tau)|d\tau\right\} - 1$$

$$\leq \exp\left\{\frac{1}{m_j}\int_{-\infty}^{t} \Omega(\tau)^{-1} |A(\tau)|d\tau\right\} - 1 \leq O(m_j^{-1}), \quad \forall t \in (-\infty, N], \tag{9.3}$$

where the constant in $O(.)$ may depend on $N$. From (9.2)

$$-\frac{4}{9}\int_{-\infty}^{a} \Omega(t)^{-2} A(t)|\psi_j(t)|^2\ dt$$

$$-\frac{4}{9}\int_{b}^{N} \Omega(t)^{-2} A(t)|\psi_j(t)|^2\ dt$$

$$\leq 4\int_{a}^{b} \Omega(t)^{-2} A(t)|\psi_j(t)|^2\ dt.$$

For a sufficiently large $j$ this implies (see (9.3))

$$-\int_{-\infty}^{a} \Omega(t)^{-2} A(t)\ e^{2m_j t} dt$$

$$-\int_b^N \Omega(t)^{-2} A(t) e^{2m_j t}\, dt$$

$$\leq \frac{45}{4} \int_a^b \Omega(t)^{-2} A(t)\, e^{2m_j t} dt.$$

This leads to a contradiction for a fixed large $N \geq \log M$. The exponential growth rate of

$$-\int_b^N \Omega(t)^{-2} A(t) e^{2m_j t}\, dt$$

as $j \to \infty$ exceeds that of the right side. By contradiction, the case $\mu_0 = 0$ is ruled out.

Next, we claim that $\mu_0 = \Omega(-\infty)$ cannot happen. Assume the contrary, $\mu_0 = \Omega(-\infty)$, i.e., there exists a sequence $\mu_j \to \Omega(-\infty)$, $\mu_j \in \mathcal{U}_{m_j}$, $m_j \to \infty$ as $j \to \infty$, $j = 1,2,3,\ldots$ . From (9.1)

$$\int_{-\infty}^{\infty} \left(|\Omega(t) - Re\mu_j|^2 + |Im\,\mu_j|^2\right)^{-1} A(t)|\psi_j(t)|^2\, dt = 0, \quad j = 1,2,3,\ldots .$$

We normalize $\psi_j$ so that

$$\psi_j(t) \sim e^{-m_j t} \text{ as } t \to \infty .$$

Selecting any $N \leq a$, we get

$$-\int_N^a \left(|\Omega(t) - Re\mu_j|^2 + |Im\,\mu_j|^2\right)^{-1} A(t)|\psi_j(t)|^2\, dt$$

$$-\int_b^{\infty} \left(|\Omega(t) - Re\mu_j|^2 + |Im\,\mu_j|^2\right)^{-1} A(t)|\psi_j(t)|^2\, dt$$

$$\leq \int_a^b \left(|\Omega(t) - Re\mu_j|^2 + |Im\,\mu_j|^2\right)^{-1} A(t)|\psi_j(t)|^2\, dt . \tag{9.4}$$

For a fixed $N \leq a$ let $j$ be large enough, so that $|\mu_j - \Omega(-\infty)| < \frac{1}{2}(\Omega(-\infty) - \Omega(N))$.

According to Lemma 2.1

$$\psi_j(t) = e^{-m_j t}\left(1 + z_j(t)\right),$$

where

$$|z_j(t)| \leq exp\left\{\frac{1}{2m_j}\int_t^{\infty} |\Omega(\tau) - \mu_j|^{-1} |A(\tau)|d\tau\right\} - 1$$

$$\leq exp\left\{\frac{1}{m_j}\int_t^\infty (\Omega(-\infty) - \Omega(\tau))^{-1} |A(\tau)|d\tau\right\} - 1 \leq O(m_j^{-1}), \qquad \forall t \in [N, \infty), \qquad (9.5)$$

where the constant in $O(.)$ may depend on $N$. Therefore, from (9.4)

$$-\frac{4}{9}\int_N^a (\Omega(-\infty) - \Omega(\tau))^{-2} A(t)|\psi_j(t)|^2\, dt$$

$$-\frac{4}{9}\int_b^\infty (\Omega(-\infty) - \Omega(\tau))^{-2} A(t)|\psi_j(t)|^2\, dt$$

$$\leq 4 \int_a^b (\Omega(-\infty) - \Omega(\tau))^{-2} A(t)|\psi_j(t)|^2\, dt\ .$$

Using (9.5) for $j$ sufficiently large, we get

$$-\int_N^a (\Omega(-\infty) - \Omega(\tau))^{-2} A(t)\, e^{-2m_j t} dt$$

$$-\int_b^\infty (\Omega(-\infty) - \Omega(\tau))^{-2} A(t) e^{-2m_j t}\, dt$$

$$\leq \frac{45}{4}\int_a^b (\Omega(-\infty) - \Omega(\tau))^{-2} A(t)\, e^{-2m_j t} dt. \qquad (9.6)$$

This leads to contradiction by selecting the limit of integration $N$ close enough to $-\infty$ and comparing the asymptotic behavior of the left side in (9.6) with the right side in (9.6) as $j \to \infty$. Thus the case $\mu_0 = \Omega(-\infty)$ is ruled out by contradiction.

The remaining possibility is $\mu_j \to \mu_0 \in (0, \Omega(-\infty))$, where $\mu_j \in \mathcal{U}_{m_j}$, $\forall j = 1,2,3,\ldots;\ m_j \to \infty$ as $j \to \infty$. For every $j$ there is an eigenfunction $\psi_j \in C^\sigma(\mathbb{R}) \cap L^2(\mathbb{R})$ exponentially decaying at $\pm\infty$. Here a fixed $\sigma \in (0,1)$. We have

$$-\frac{d^2}{dt^2}\psi_j + m_j^2 \psi_j + (\Omega(t) - \mu_j)^{-1} A(t)\psi_j = 0, \ \ t \in \mathbb{R}\ ,\ Im\, \mu_j > 0,\ \psi_j \neq 0, \qquad (9.7)$$

We normalize $\psi_j$ so that

$$\|\psi_j\|_{C^\sigma(\mathbb{R})} = 1\ .$$

Splitting $m_j^2 = (m_j^2 - m_0^2) + m_0^2$ for some $m_0 > 0$ we get in $C^\sigma(\mathbb{R})$

$$\psi_j + (m_j^2 - m_0^2)K_{m_0}\psi_j + K_{m_0}(\Omega(t) - \mu_j)^{-1} A(t)\psi_j = 0. \qquad (9.8)$$

We claim the sequence

$$\{j \to K_{m_0}(\Omega(t) - \mu_j)^{-1} A(t)\psi_j\}$$

is relatively compact in $C^{\sigma_2}(\mathbb{R})$ for some $\sigma_2 \in (\sigma, 1)$. Indeed, for any $p \in [1, \infty]$,

$\|A(t)\psi_j\|_{L^p(\mathbb{R})} \leq C < \infty$. Let as before $\mu_0 = \Omega(v_0), v_0 \in \mathbb{R}$. We select $\chi \in C_0^\infty(\mathbb{R})$ such that $\chi(t) = 1$ for $t \in [v_0 - 1, v_0 + 1]$ and split

$$I_j(t) = K_{m_0}(\Omega(t) - \mu_j)^{-1} A(t)\psi_j$$

$$= K_{m_0}\chi(t)(\Omega(t) - \mu_j)^{-1} A(t)\psi_j + K_{m_0}(1 - \chi(t))(\Omega(t) - \mu_j)^{-1} A(t)\psi_j$$

$$= I_{1j}(t) + I_{2j}(t).$$

We have $\|I_{2j}\|_{C^{\sigma_1}(\mathbb{R})} \leq C < \infty$ for a fixed $\sigma_1 \in (\sigma, 1)$; $\|I_{2j}\|_{L^p(\mathbb{R})} \leq C < \infty$ for any $p \in (1, \infty)$.

Moreover, from the explicit form of the kernel $K_{m_0}$ and a pointwise inequality $|A(t)\psi_j(t)| \leq |A(t)|$, $\forall t \in \mathbb{R}$, we derive that

$$\forall \varepsilon > 0 \; \exists R > 0 \text{ such that } \|I_{2j}\|_{L^p(\mathbb{R} \setminus [-R,R])} < \varepsilon, \; j = 1,2,3,\ldots \; .$$

Indeed, pointwise

$$|I_{2j}(t)| \leq C \frac{1}{2m_0} \int e^{-m_0|t-\xi|} |A(\xi)| d\xi$$

and the statement follows. We turn attention to the first term $I_{1j}(t)$.

Examining the proof of Lemma 6.1 we find

$$\|I_{1j}\|_{C^{\sigma_1}(\mathbb{R})} \leq C < \infty, \; j = 1,2,3,\ldots ,$$

and

$$\forall \varepsilon > 0 \; \exists R > 0 \text{ such that } \|I_{1j}\|_{L^p(\mathbb{R} \setminus [-R,R])} < \varepsilon, \; j = 1,2,3,\ldots \; .$$

Combining the estimates for $I_{1j}$ and $I_{2j}$ yields

$$\|I_j\|_{C^{\sigma_1}(\mathbb{R})} \leq C < \infty, \; j = 1,2,3,\ldots ,$$

and the sequence $\{I_j\}, j = 1,2,3,\ldots$ is relatively compact in $L^p(\mathbb{R})$. therefore, after possibly extracting a subsequence we may assume that $\{I_j\}, j = 1,2,3,\ldots$ converges strongly in $L^p(\mathbb{R})$ and the limit belongs to $C^{\sigma_1}(\mathbb{R})$. Therefore, this convergence is strong in $C^{\sigma_2}(\mathbb{R})$ for any $\sigma_2 \in [\sigma, \sigma_1)$.

Therefore, after possibly extracting a subsequence we may assume that

$$K_{m_0}\big(\Omega(t) - \mu_j\big)^{-1} A(t)\psi_j \to f \text{ strongly in } C^{\sigma_2}(\mathbb{R}) \text{ and in } L^p(\mathbb{R})$$

for some $\sigma_2 \in [\sigma, 1)$.

We need the following lemma.

Lemma 9.1. For any fixed $s \in \mathbb{R}$ the pseudo-differential operator

$$f \to Bf \equiv (2\pi)^{-1} \int\int (\xi^2 + m^2)^{-1}(\xi^2 + m_0^2) e^{i\xi(x-y)} f(y)\, dy d\xi$$

Is bounded in $C^s(\mathbb{R})$ uniformly in $m \in [1, \infty)$.

Remark. For an integer $s$, $C^s(\mathbb{R})$ is understood in the sense of Zygmund.

Proof. Let

$$f = \sum_{l=-1}^{\infty} \Delta_l f$$

be Littlewood-Paley decomposition of $f \in C^s(\mathbb{R})$.

We estimate

$$\kappa_l(z) = (2\pi)^{-1} \int \varphi(2^{-l}\xi)(\xi^2 + m^2)^{-1}(\xi^2 + m_0^2) e^{i\xi z}\, d\xi.$$

We have for $k = 0,1,2,3,\dots;\ l = 0,1,2,3,\dots$ .

$$(-iz)^k \kappa_l(z) = (2\pi)^{-1} \int \partial_\xi^k \big(\varphi(2^{-l}\xi)(\xi^2 + m^2)^{-1}(\xi^2 + m_0^2)\big) e^{i\xi z}\, d\xi.$$

Therefore,

$$|z^k \kappa_l(z)| \le C_k 2^{l(1-k)}, \quad \forall z \in \mathbb{R},$$

where the constant $C_k$ does not depend on $l \ge 0, m \ge 1$. Thus,

$$\|\kappa_l\|_{L^1(\mathbb{R})} \le C \int_{|z| \le 2^{-l}} 2^l\, dz + C \int_{|z| \ge 2^{-l}} |z|^{-2} 2^{-l}\, dz \le C,$$

where the constant $C$ does not depend on $l \ge 0, m \ge 1$. We have

$$\Delta_l Bf = \kappa_l * f = \sum_{|j-l| \le 1} \kappa_l * \Delta_j f.$$

Therefore,

$$\|\Delta_l Bf\|_{L^\infty(\mathbb{R})} \le C \sum_{|j-l| \le 1} \|\Delta_j f\|_{L^\infty(\mathbb{R})} \le C \|f\|_{C^s(\mathbb{R})} 2^{-sl}, \quad l = 0,1,2,3,\dots\ .$$

Here the constant $C$ does not depend on $l \geq 0, m \geq 1$.

For $l = -1$, using Bernstein's inequality, we get

$$\|\Delta_l Bf\|_{L^\infty(\mathbb{R})} = \|K_m(-\tfrac{d^2}{dt^2}+m_0^2)\Delta_{-1}f\|_{L^\infty(\mathbb{R})} \leq Cm^{-2}\|\Delta_{-1}f\|_{L^\infty(\mathbb{R})}$$

$$\leq Cm^{-2}\|f\|_{C^s(\mathbb{R})}.$$

This completes the proof of Lemma 9.1. QED.

Using (9.8) for $l = -1,0,1,2,\ldots$

$$\Delta_l \psi_j = -B\Delta_l K_{m_0}(\Omega(t) - \mu_j)^{-1} A(t)\psi_j,$$

where $m = m_j$. We choose a large $N$ and using boundedness of the right side in $C^{\sigma_1}(\mathbb{R})$ we get the inequality

$$\|\Delta_l \psi_j\|_{L^\infty(\mathbb{R})} \leq C 2^{-l\sigma_1}, \quad \forall j = 1,2,3,\ldots; \forall l \geq N. \tag{9.9}$$

On the other hand,

$$\Delta_l \psi_j = -K_{m_j}\left(-\frac{d^2}{dt^2} + m_0^2\right)\Delta_l K_{m_0}(\Omega(t) - \mu_j)^{-1} A(t)\psi_j.$$

Using Bernstein's inequality and the identity

$$\|(2m_j)^{-1} e^{-m_j|\cdot|}\|_{L^1(\mathbb{R})} = m_j^{-2}, \quad j = 1,2,3,\ldots,$$

we get

$$\|\Delta_l \psi_j\|_{L^\infty(\mathbb{R})} \leq Cm_j^{-2} 2^{2l} 2^{-l\sigma_1}, \quad \forall j = 1,2,3,\ldots; \forall l = -1,0,1,\ldots,N. \tag{9.10}$$

From (9.9), (9.10) for sufficiently large $j$ we get

$$2^{l\sigma}\|\Delta_l \psi_j\|_{L^\infty(\mathbb{R})} \leq C 2^{-N(\sigma_1 - \sigma)}, \forall j \gg 1; \forall l = -1,0,1,\ldots.$$

Therefore,

$$\|\psi_j\|_{C^\sigma(\mathbb{R})} \leq C 2^{-N(\sigma_1 - \sigma)}, \quad \forall j \gg 1.$$

This leads to a contradiction for $N$ large enough with normalization $\|\psi_j\|_{C^\sigma(\mathbb{R})} = 1$ thus completing the proof of Proposition 9.1. QED.

## §9. Structure of the spectrum

We study the family of operators

$$L_\mu + m^2 \equiv -\frac{d^2}{dt^2} + (\Omega(t) - \mu)^{-1} A(t) + m^2 \tag{10.1}$$

as $m$ decreases from $\infty$ to $1$.

Theorem 10.1. Let $-m_d^2$, $m_d > 1$, $d = a$ or $d = b$, be the bottom of the spectrum in $L^2(\mathbb{R})$ of

$$L_{\Omega(d)} = -\frac{d^2}{dt^2} + (\Omega(t) - \Omega(d))^{-1} A(t). \tag{10.2}$$

Then,

$$\mathcal{U}_m = \emptyset, \quad \forall m \in [m_a, \infty);$$

$$\#\mathcal{U}_m = 1, \quad \forall m \in (m_b, m_a);$$

$$\mathcal{U}_m = \emptyset, \quad \forall m \in [1, m_b].$$

For the corresponding eigenfunction (unique up to a nonzero complex factor) $\psi \in L^2(\mathbb{R})$

$$L_\mu \psi + m^2 \psi = 0, \quad m \in (m_b, m_a), \quad \mathcal{U}_m = \{\mu\},$$

we have

$$\int (\Omega(t) - \mu)^{-2} A(t) \psi^2(t)\, dt \neq 0. \tag{10.3}$$

Remark 10.1. See Propositions 4.1 and 4.2.

The proof of Theorem 10.1 relies on the information obtained above and on the following statement.

Proposition 10.1. (Stability of the root multiplicities). Let $A \in \mathcal{L}(H)$ be a bounded operator in a Hilbert space $H$, and let $\Gamma$ be a rectifiable closed curve with interior $D$, with positive orientation. Assume that

1. Every point of $\Gamma$ is a regular point of $A$;
2. $\sigma_A^{ess} \cap D = \emptyset$.

Define

$$P_\Gamma = \frac{1}{2\pi i} \oint_\Gamma (z - A)^{-1} dz\,;$$

$$\nu_\Gamma = \nu_\Gamma(A) = \dim P_\Gamma H < \infty\,.$$

Then, there exists a $\delta > 0$, such that for any $B \in \mathcal{L}(H)$, $\|A - B\|_{\mathcal{L}(H)} < \delta$, both conditions 1. and 2. hold, and

$$\nu_\Gamma(B) = \nu_\Gamma(A)\,.$$

We make one general remark concerning the eigenfunction $\psi \in L^2(\mathbb{R})$ of the operator $L_\mu$

$$L_\mu \psi + m^2 \psi = 0, \quad m > 0,\ Im\,\mu > 0. \tag{10.4}$$

It follows from Lemma 2.1 that $\psi$ decays exponentially at $\pm\infty$ and belongs to a Holder space $C^\sigma(\mathbb{R})$, $\forall \sigma \in \mathbb{R}$. It follows from (10.1), (10.4) that for

$$g = -\frac{d^2}{dt^2}\psi + m^2 \psi, \tag{10.5}$$

we get

$$(\Lambda_m g)(t) \equiv \Omega(t) g(t) + A(t)(K_m g)(t) = \mu g(t). \tag{10.6}$$

The operator $\Omega(t) + A(t) K_m$ in the left side of (10.6) is bounded in $L^2(\mathbb{R})$, and every solution

$$g \in L^2(\mathbb{R})$$

of the equation (10.6) with $Im\,\mu > 0$ will produce a solution $\psi = K_m g \in L^2(\mathbb{R})$ to the equation (10.4).

Proof of the Theorem 10.1. From Proposition 9.1 above, $\mathcal{U}_m = \emptyset$ for $m \geq m_*$, where $m_*$ is sufficiently large. It follows from (10.6) and from identity, $\|K_m\|_{\mathcal{L}(L^2(\mathbb{R}))} = m^{-2}$ that every $\mu \in \mathcal{U}_m$ satisfies

$$|\mu| \leq R,\ m \geq \frac{1}{2}.$$

We claim $\mathcal{U}_m = \emptyset$ for $m \in [m_a, \infty)$. Indeed, assume $m \in [m_a, \infty)$ and $\exists \mu \in \mathcal{U}_m$. Using Proposition 10.1 we construct a sequence $m_j \to m_0 \in (m_a, m_*]$ and a sequence $\mu_j \in \mathcal{U}_{m_j}$ such that $\mu_j \to \mu_0 \in \mathbb{R}$, as $j \to \infty; j = 1,2,3, \ldots$ .

If $\mu_0 \in \mathbb{R}\setminus[0,\Omega(-\infty)]$ this contradicts Proposition 3.4. The case $\mu_0 = 0$ is ruled out by Proposition 3.1. The case $\mu_0 = \Omega(-\infty)$ cannot happen according to Proposition 3.2. Therefore, $\mu_0 \in (0, \Omega(-\infty))$. From Proposition 3.3 this implies $\mu_0 = \Omega(d)$, where either $d = a$ or $d = b$.

From propositions 4.1 and 4.2 we have

$$m_0 = m_a \text{ or } m_0 = m_b, \tag{10.7}$$

depending on whether $d = a$ or $d = b$.

But by construction $m_0 > m_a$. Also, $m_a > m_b$ thus $m_0 > m_b$. This contradicts (10.7). Therefore,

$$\mathcal{U}_m = \emptyset \, , \, \forall m \in [m_a, \infty) . \tag{10.8}$$

We claim

$$\lim_{m \to m_a - 0} \#\mathcal{U}_m = 1. \tag{10.9}$$

Indeed, we constructed in the proof of Theorem 5.1 an eigenfunction $\psi_{m,\mu}$ with $m = m_0 - h, h \in (0, \varepsilon)$ and $|\mu - \mu_0| < Ch$, $\text{Im}\,\mu > 0$. Here $\mu_0 = \Omega(a)$. This is a unique solution in the ball $\{(|\mu - \mu_0|^2 + |m - m_0|^2)^{\frac{1}{2}} < \delta\}$ for a sufficiently small $\delta > 0$. Indeed, let there exists for every $\delta = i^{-1}$, positive integer $i \geq N$ sufficiently large, a solution $(m_i, \mu_i') \to (m_0, \mu_0)$ different from the one constructed there, i.e., $\mu_i' \neq \mu_i$; $\text{Im}\,\mu_i' > 0$; $(|\mu_i' - \mu_0|^2 + |m_i - m_0|^2)^{\frac{1}{2}} < i^{-1}$. This means

$$(\psi_{m_i,\mu_i}, \psi_0) - 1 = 0; \, (\psi_{m_i,\mu_i'}, \psi_0) - 1 = 0, \quad i = N, N+1, N+2, \ldots .$$

We have the identity (7.2) valid for both solutions $(m_i, \mu_i)$ and $(m_i, \mu_i')$. Subtracting we get

$$(\mu_i' - \mu_i)((\Omega - \mu_0 - i0)^{-1}[(\Omega - \mu_0)^{-1}A]\psi_0, \psi_0)$$

$$= o((|m_i - m_0|^2 + |\mu_i - \mu_0|^2 + |\mu_i' - \mu_0|^2)^{\frac{1}{2}})$$

$$= o\left((|m_i - m_0|^2 + |\mu_i' - \mu_i|^2)^{\frac{1}{2}}\right), \tag{10.9}$$

since

$$|\mu_i' - \mu_0| \leq |\mu_i' - \mu_i| + |\mu_i - \mu_0|; \, |\mu_i - \mu_0| = O(|m_i - m_0|), \quad i = N, N+1, N+2, \ldots .$$

The constant in the left side of (10.9) is not zero, thus

$$|\mu_i' - \mu_i| = o(|m_i - m_0|), \quad i = N, N+1, N+2, \ldots .$$

This statement puts $\mu_i'$ inside the circle, where Rouchet's theorem there is one and only one solution to the equation

$$(\psi_{m_i,\mu}, \psi_0) - 1 = 0 ,$$

namely $\mu = \mu_i$. Therefore,

$$\mu_i' = \mu_i , \qquad , \qquad i = N, N+1, N+2, \dots .$$

The assumption $(m_i, \mu_i') \to (m_0, \mu_0)$ is a valid one because of the argument above with deformation of the spectrum which would lead to

$$\mathcal{U}_{m_a} \neq \emptyset$$

otherwise, in contradiction with (10.8). From (10.9) it follows that

$$\# \mathcal{U}_m = 1 , \quad \forall\, m \in (m_b, m_a). \tag{10.10}$$

The branch of solutions with $m = m_b + h, h \in (0, \varepsilon)$ and $|\mu - \mu_0| < Ch$, $Im\,\mu > 0$ has been constructed in the proof of Theorem 5.1. Here $\mu_0 = \Omega(b)$. In combination with (10.10) the deformation argument yields

$$\mathcal{U}_{m_b} = \emptyset . \tag{10.11}$$

Applying again the deformation argument to the interval $m \in [1, m_b)$ and using the study of the branch near the point $(m_b, \Omega(b))$ and the statement (10.11) we conclude that

$$\mathcal{U}_m = \emptyset , \qquad \forall m \in [1, m_b). \tag{10.12}$$

We constructed a unique (up to a nonzero complex factor) eigenfunction $\psi \in L^2(\mathbb{R})$ of the problem

$$L_\mu \psi + m^2 \psi = 0 , \qquad m \in (m_b, m_a) , \qquad \mathcal{U}_m = \{\mu\}.$$

By construction and by the deformation of spectrum argument this solution is of algebraic multiplicity 1. For the spectral problem (10.6), (10.5), where $g \in L^2(\mathbb{R})$, $\psi = K_m g$ we conclude that the equation

$$(L_\mu + m^2)\psi_1 = (\Omega(t) - \mu)^{-1}(-\frac{d^2}{dt^2} + m^2)\psi$$

Has no solution $\psi_1 \in L^2(\mathbb{R})$. But

$$(-\frac{d^2}{dt^2} + m^2)\psi = -(\Omega(t) - \mu)^{-1} A(t)\psi$$

Also, the null space of the adjoint operator $(L_\mu + m^2)^*$ is spanned by $\bar\psi$.

This implies

$$\int (\Omega(t) - \mu)^{-1} A(t) \psi(t)^2 dt \neq 0 , \qquad (10.13)$$

as claimed. Combining (10.8), (10.10), (10.11), (10.12), (10.13) we arrive at the statement of Theorem 10.1. QED.

Remark. A simple but lengthy argument shows that the function $(m_b, m_a) \to H$, $\{m \to \mu(m)\}$, where

$$\mathcal{U}_m = \{\mu(m)\}$$

is differentiable with

$$\frac{d\mu(m)}{dm} = -2m \left( \int \psi(t)^2 dt \right) \left( \int (\Omega(t) - \mu)^{-1} A(t) \psi(t)^2 dt \right)^{-1} ,$$

see (10.13).

# §11. Construction of the velocity profile

Theorem 11.1. For any $\alpha \in (0,2)$ there exists a function $\Omega \in \mathcal{C}$ and an *integer* $m \geq 2$ so that

$$\#\mathcal{U}_m = 1;$$

$$\mathcal{U}_{m+l} = \emptyset, \quad \forall l = 1,2,3,\dots .$$

Proof. Let $\Omega_0 \in \mathcal{C}$ be a function such that $a = 0$, $b = 1$. We denote the bottom of the spectrum of the operator

$$L_0 \equiv -\frac{d^2}{dt^2} + (\Omega_0(t) - \Omega_0(a))^{-1} (\Omega_0''(t) + 2\Omega_0'(t)) \geq -N_0^2 \qquad (11.1)$$

in $L^2(\mathbb{R})$ by $-N_0^2$. Here we choose $N_0 > 0$.

Choose $m \geq 2$ to be an integer such that

$$m^2 > N_0^2. \tag{11.2}$$

Construct using Proposition 8.1 a function $\Omega_1 \in \mathcal{C}$ such that $a = 0, \ b = 1$ and

$$L_1 \equiv -\frac{d^2}{dt^2} + \left(\Omega_1(t) - \Omega_1(a)\right)^{-1}\left(\Omega_1''(t) + 2\Omega_1'(t)\right) \geq -N_1^2, \tag{11.3}$$

Where $-N_1^2$ ($N_1 > 0$) the bottom of the spectrum of $L_1$ in $L^2(\mathbb{R})$ satisfies the inequality

$$N_1^2 > m^2. \tag{11.4}$$

For any $\theta \in [0,1]$ let

$$\Omega_\theta = \theta\Omega_1 + (1-\theta)\Omega_0 \in \mathcal{C}.$$

We denote as $-N_\theta^2$ ($N_\theta > 0$) the bottom of the spectrum in $L^2(\mathbb{R})$ of the operator

$$L_\theta \equiv -\frac{d^2}{dt^2} + \left(\Omega_\theta(t) - \Omega_\theta(a)\right)^{-1}\left(\Omega_\theta''(t) + 2\Omega_\theta'(t)\right) \geq -N_\theta^2. \tag{11.5}$$

Also, let

$$M_\theta \equiv -\frac{d^2}{dt^2} + \left(\Omega_\theta(t) - \Omega_\theta(b)\right)^{-1}\left(\Omega_\theta''(t) + 2\Omega_\theta'(t)\right) \geq -W_\theta^2, \tag{11.6}$$

where $-W_\theta^2$ ($W_\theta > 0$) stands for the bottom of the spectrum in $L^2(\mathbb{R})$ of the operator $M_\theta$. The functions

$$\theta \to N_\theta$$

$$\theta \to W_\theta$$

are continuous on $[0,1]$. Also,

$$W_\theta < N_\theta, \forall \theta \in [0,1]. \tag{11.7}$$

Let

$$\theta_0 = \inf\{\theta \in [0,1] | \ N_\theta > m \text{ on } (\theta, 1]\}. \tag{11.8}$$

From (11.1) --(11.8) $N_{\theta_0} = m$ and $N_\theta > m$, $\forall \theta \in (\theta_0, 1]$. Since $W_{\theta_0} < m$, there is a $\delta \in (0, 1 - \theta_0)$ such that $W_{\theta_0+\delta} < m$, $N_{\theta_0+\delta} < m + 1$.

Therefore,

$$W_{\theta_0+\delta} < m < N_{\theta_0+\delta} < m + 1.$$

For $\Omega = \Omega_{\theta_0+\delta}$ we have from Theorem 10.1

$$\#\boldsymbol{U}_m = 1;$$

$$\boldsymbol{U}_{m+l} = \emptyset, \quad \forall l = 1,2,3,\ldots\ .$$

This concludes the proof of Theorem 11.1.   QED.